\documentclass[a4paper]{article}

\usepackage[english]{babel}
\usepackage[utf8x]{inputenc}
\usepackage[T1]{fontenc}
\usepackage[a4paper,top=3cm,bottom=2cm,left=3cm,right=3cm,marginparwidth=1.75cm]{geometry}
\usepackage{amsmath}
\usepackage{graphicx}
\usepackage[colorinlistoftodos]{todonotes}
\usepackage{listings}
\usepackage{color} 
\definecolor{mygreen}{RGB}{28,172,0} 
\definecolor{mylilas}{RGB}{170,55,241}
\usepackage{changepage}
\usepackage{graphicx}
\usepackage{caption}
\usepackage{subcaption}
\usepackage[makeroom]{cancel}
\usepackage{listings}
\usepackage{amsmath}
\usepackage{amsfonts}
\usepackage{float}
\usepackage{listings}
\usepackage{framed}
\usepackage{cite}
\usepackage{algorithm}
\usepackage{pdfpages}
\usepackage{mathtools}
\usepackage{bm}
\usepackage{authblk}

\usepackage{todonotes}

\newcommand{\norm}[1]{\left\lVert#1\right\rVert}

\providecommand{\keywords}[1]{\textit{Keywords and phrases: } #1}

\begin{document}
	
\title{Efficiency of a micro-macro acceleration method for scale-separated stochastic differential equations}

\author[1]{Hannes Vandecasteele}
\author[1]{Przemys{\l}aw Zieli{\'n}ski}
\author[1]{Giovanni Samaey}
\affil[1]{KU Leuven, Department of Computer Science, NUMA Section, Celestijnenlaan 200A box 2402, 3001 Leuven, Belgium}
\date{\today}
	
\maketitle

\begin{abstract}
We discuss through multiple numerical examples the accuracy and efficiency of a micro-macro acceleration method for stiff stochastic differential equations (SDEs) with a time-scale separation between the fast microscopic dynamics and the evolution of some slow macroscopic state variables. The algorithm interleaves a short simulation of the stiff SDE with extrapolation of the macroscopic state variables over a longer time interval. After extrapolation, we obtain the reconstructed microscopic state via a matching procedure: we compute the probability distribution that is consistent with the extrapolated state variables, while minimally altering the microscopic distribution that was available just before the extrapolation. In this work, we numerically study the accuracy and efficiency of micro-macro acceleration as a function of the extrapolation time step and as a function of the chosen macroscopic state variables. Additionally, we compare the effect of different hierarchies of macroscopic state variables. We illustrate that the method can take significantly larger time steps than the inner microscopic integrator, while simultaneously being more accurate than approximate macroscopic models. 
\end{abstract}


\section{Introduction}
Many applications in science and engineering are modeled with stochastic differential equations (SDEs)
\begin{equation} \label{eq:SDE}
d\mathbf{X}_t = a(\mathbf{X}_t,t)dt + b(\mathbf{X}_t,t)d\mathbf{W}_t,
\end{equation}
where $\mathbf{X}_t$ is a diffusion process on a domain $G \subset \mathbb{R}^d$. The drift term $a(\bm{X}_t, t) \in \mathbb{R}^d$ has the same dimension as $\bm{X}_t$, the diffusion term $b(\bm{X}_t,t) \in \mathbb{R}^{d \times n}$ and $\mathbf{W}_t \in \mathbb{R}^n$ represents an $n$-dimensional Brownian motion.
In practice however, we are often only concerned with the evolution of a few macroscopic state variables $\bm{m}(t) = \left(m_1(t), \dots, m_L(t)\right)$, which are defined as the expectation of some functions of interest $\bm{R}(\bm{x}) = \left(R_1(\bm{x}), \dots, R_L(\bm{x})\right)$ over the distribution of $\mathbf{X}_t$, i.e.,
\begin{equation} \label{eq:stateevolution}
t \mapsto \bm{m}(t) = \mathbb{E}[\bm{R}(\mathbf{X}_t)].
\end{equation}
 In many situations, it is difficult or impossible to derive a closed model for the dynamics of the macroscopic states $\bm{m}(t)$. Instead, we need to approximate the expectation in~\eqref{eq:stateevolution} at different points in time via Monte Carlo methods and simulate the SDE~\eqref{eq:SDE} for each of the generated samples. 

 In this work, we are concerned with SDE systems that are stiff, in the sense that individual microscopic paths vary over fast time scales, while the evolution of the macroscopic state variables takes place over much longer time scales. Due to the discrepancy in time scales, explicit time discretization methods, such as the Euler-Maruyama scheme, are forced to take small time steps as they have a restricted stability domain. Hence, explicit methods are inadequate to simulate the macroscopic state variables up to a large end time $T$.

Over the years, many multiscale methods have been proposed to overcome the problem of stiffness. We mention implicit methods \cite{tian2001implicit,amiri2015class}, S-ROCK \cite{abdulle2008s}, the equation-free framework \cite{kevrekidis2009equation,kevrekidis2003equation} and the heterogeneous multiscale method (HMM) \cite{weinan2003multiscale,abdulle2012heterogeneous}. Implicit methods have proven to be very successful for stiff ODEs since they allow for much larger time steps than their explicit counterparts, due to a much larger stability domain, \cite{byrne1987stiff}. The cost of implicit time steppers is usually higher, since at every time step a (non-)linear system needs to be solved. Nevertheless, implicit schemes are a major improvement over explicit time integrators on stiff ODEs because they allow taking much larger time steps. However, a problem arises when applying implicit techniques to SDEs \cite{burrage2004numerical}. The authors of \cite{li2008effectiveness} showed that these methods are unable to capture the probability distribution of the fast modes when taking large time steps. The idea of S-ROCK is to extend the stability domain maximally along the negative real axis by bounding the stability domain by a Chebyshev polynomial \cite{abdulle2008s,abdulle2007stabilized,komori2012weak}. One can then construct a Runge-Kutta scheme that has this stability domain. A drawback is that S-ROCK schemes do not attain a high convergence order \cite{komori2013strong}.

Both the equation-free framework \cite{kevrekidis2009equation} and the heterogeneous multiscale method \cite{abdulle2012heterogeneous} can yield considerable simulation speed-ups. However, the convergence analysis is usually restricted to showing converge in the limit when the time-scale separation becomes infinite. In that limit, the microscopic dynamics itself converges to a limiting macroscopic model, and both the equation-free method and HMM recover that limiting macroscopic model when the time-scale separation becomes infinite\cite{abdulle2012heterogeneous}. An alternative approach is to average out the fast dynamics, either analytically or numerically, to obtain a the evolution of the slow dynamics only \cite{legoll2010effective,pavliotis2008multiscale}. All of these methods introduce a modeling error when the time-scale separation is finite. However, in cases of large time-scale separation, the limiting macroscopic model becomes accurate and this modeling error becomes small, compared to the time discretization error.  

Recently, a new micro-macro acceleration algorithm \cite{debrabant2017micro} was introduced as an alternative to the above techniques. Micro-macro acceleration exploits the stiffness of the SDE by introducing a second time step, besides the microscopic time step $\delta t$, which is used to simulate the stiff SDE~\eqref{eq:SDE}. This second, much larger time step $\Delta t \gg \delta t$, is introduced to extrapolate the evolution of the macroscopic state variables $(m_1(t),\dots,m_L(t))$. One time step of micro-macro acceleration consists of four steps: (i) Monte Carlo \textit{simulation} of the SDE~\eqref{eq:SDE} with a small time step $\delta t$; (ii) \textit{restriction} to approximate the macroscopic state variables~\eqref{eq:stateevolution} at every microscopic time step; (iii) \textit{extrapolation} of the macroscopic state variables over a time interval of size $\Delta t$; and (iv) \textit{matching} to generate a new microscopic probability distribution, consistent with the extrapolated macroscopic states, and with minimal deviation from a \textit{prior} distribution (which we choose to be the final microscopic distribution obtained during the simulation stage). There are many ways of expressing the deviation of one probability distribution from another. One possibility is minimizing the $L_2$-norm between the matched and prior distribution, but this choice does not guarantee positivity of the matched distribution \cite{debrabant2017micro}. In this text, we consider matching in Kullback-Leibler divergence, which is based on notions from information theory \cite{lelievre2018analysis}. We provide a more detailed mathematical description of the algorithm in Section 2.

It is shown that micro-macro acceleration can converge to the exact dynamics of the SDE, even when fixing a finite time-scale separation, for both matching in the $L_2$-norm \cite{debrabant2017micro} and in Kullback-Leibler divergence \cite{lelievre2018analysis}. The conditions under which convergence has been demonstrated are the following: (i) the number of state variables $L$ must increase to infinity; and (ii) the time steps present (the microscopic time step $\delta t$ and the extrapolation time step $\Delta t$), must tend to zero \cite{lelievre2018analysis,debrabant2017micro}. Additionally, it was shown that the method is stable for the extrapolation step sizes independent of the fast time scales in the system \cite{debrabant2018study}. 

The convergence analysis in \cite{debrabant2017micro,lelievre2018analysis} reveals that the errors in the micro-macro acceleration method that are caused by the combination of extrapolation and matching, can be viewed as a time discretization error, with a size that depends on the extrapolation time step. In contrast, equation-free \cite{kevrekidis2009equation,kevrekidis2003equation} and heterogeneous multiscale methods \cite{weinan2003multiscale,abdulle2012heterogeneous} introduce a modeling error, because they compute the time evolution of an approximate macroscopic model in which a numerical closure relation is introduced. Since micro-macro acceleration is stable for time steps that do not depend on the fast time scales, the extrapolation time step can be chosen based on accuracy considerations only. Hence, we need to study accuracy as a function of extrapolation time step to assess the efficiency of the method. This study is exactly the aim of the current work. 

The main parameters that determine the accuracy of micro-macro acceleration are the following:
\begin{itemize}
\item The number of macroscopic state variables $L$. We may expect that for a larger value of $L$, the matched distribution lies closer to the exact microscopic distribution, since the method converges as $L$ tends to infinity. However, matching is computationally more expensive when $L$ increases. 
\item The choice of the functions of interest $R_l$ in equation~\eqref{eq:stateevolution} and their potential to capture the essential features of the underlying microscopic distribution. Two different hierarchies of macroscopic state variables may reach a different accuracy using the same number $L$ of state variables, due to the different nature of both hierarchies.
\item The extrapolation step size $\Delta t$. When $\Delta t$ decreases, we may expect that micro-macro acceleration is more accurate, but the computational cost increases.
\end{itemize}

Based on the arguments above, we choose to study the accuracy that can be reached for a given number of state variables $L$ and extrapolation step size $\Delta t$. We investigate the choice of $L$ and $\Delta t$ on a set of representative toy examples. In particular, we investigate for which parameters $L$ and $\Delta t$ micro-macro acceleration can be more accurate than the approximate macroscopic model (or numerical closure), while simultaneously being faster than a full microscopic simulation.

The text is organized as follows. In Section~\ref{sec:mM}, we describe the micro-macro acceleration algorithm in mathematical detail. We discuss some implementation issues in Section~\ref{sec:impl_det}. Section~\ref{sec:states} is devoted to the choice of the hierarchy and the number of macroscopic state variables to extrapolate. We focus on two problems from molecular dynamics: FENE dumbbells and a simple three-atom molecule. In Section~\ref{sec:extraoplation}, we study the accuracy of micro-macro acceleration as a function of the extrapolation step $\Delta t$, when fixing the number of macroscopic state variables. More specifically, we discuss two stochastic models with a time-scale separation for which an approximate macroscopic model can be derived for the slow component of the system, in the limit of infinite time-scale separation. The first example is an overdamped Langevin equation with a double-well potential; the second is a linear system with an external time-periodic force. We investigate for which extrapolation step size micro-macro acceleration is more accurate than the corresponding approximate macroscopic models. Finally, Section~\ref{sec:conclusion} presents a concluding discussion.

\section{A micro-macro acceleration algorithm} \label{sec:mM}
In this Section, we introduce the micro-macro acceleration algorithm and delineate each of its four steps. We define the microscopic time-stepper that integrates the stochastic process \eqref{eq:SDE} over small time steps $\delta t$ using a weighted Monte Carlo ensemble in Section~\ref{subsec:mc}. In Section~\ref{subsec:restriction}, we present the restriction operator to retrieve the macroscopic state variables from the microscopic Monte Carlo particle ensemble. Then, in Section~\ref{subsec:extrapolation}, we introduce the extrapolation of the macroscopic state variables over a larger time step $\Delta t$. Finally, Section~\ref{subsec:matching} contains the description of the matching operator with which we return from the macroscopic to the microscopic level of description. Matching is performed by reweighting the Monte Carlo particle ensemble. We collect the complete algorithm in Section~\ref{subsec:algorithm}.

\subsection{Monte Carlo simulation} \label{subsec:mc}
Suppose at time $t^n = n\Delta t$ we have a weighted particle ensemble $\mathcal{X}_J^n = \left(w_j^n, \bm{X}_j^n\right)_{j=1}^J$, where the state variable $\bm{X}_j^n$ denotes the $j$-th realization of the SDE \eqref{eq:SDE} and $w_j^n$ is its associated weight. The particle ensemble $\mathcal{X}^n$ determines an empirical probability distribution $\hat{\mu}^n$
\[
\hat{\mu}^n = \Sigma_{j=1}^J w_j^n \delta_{\bm{X}_j^n},
\]
that approximates the exact continuous probability distribution $\mu^n$ of the diffusion process \eqref{eq:SDE} at time $t^n$.
In the first stage of micro-macro acceleration, we perform $K$ microscopic time steps of size $\delta t$ with a Monte Carlo time integrator. The propagation of the particle ensemble through the SDE \eqref{eq:SDE} at time $t^n+k\delta t$ is given by a (possibly time-dependent) transition operator $S^{n,k}_{\delta t}$
\begin{equation} \label{eq:microtimestepper}
\mathcal{X}_J^{n,k} = S^{n,k}_{\delta t} (\mathcal{X}_J^{n,k-1}), \ \ \ k = 1,\dots,K,
\end{equation}
where every particle ensemble $\mathcal{X}_J^{n,k}$ determines the empirical probability distribution $\hat{\mu}^{n,k} = \Sigma_{j=1}^J w_j^n \delta_{\bm{X}_j^{n,k}}$ that approximates the exact probability distribution $\mu^{n,k}$ at time $t^n+k\delta t$. For completeness, we define $\mathcal{X}_J^{n,0}=\mathcal{X}_J^{n}$. For instance, the Euler-Maruyama scheme propagates each particle $\mathbf{X}_j^{n,k} \in \mathcal{X}^{n,k}$ as
\[
\bm{X}_j^{n,k+1} = \bm{X}_j^{n,k}+a(\bm{X}_j^{n,k}, t^n+k\delta t)\delta t + \sqrt{\delta t} b(\bm{X}_j^{n,k}, t^n+k\delta t) \bm{\xi}^{n,k},
\]
with $\bm{\xi}^{n,k}$ an $n$-dimensional standard normally distributed random variable. For the remainder of the manuscript, we will use the Euler-Maruayama method as the inner time integrator.

\subsection{Restriction} \label{subsec:restriction}
To transit from the microscopic to the macroscopic description of the diffusion process, we introduce the restriction operator $\mathcal{R}$ that acts on a probability distribution $\mu$ as follows
\begin{equation} \label{eq:restrictionoperator}
\mathcal{R}(\mu) = \mathbb{E}_{\mu}[\mathbf{R}] = \int_G \mathbf{R}(\bm{x}) d\mu(\bm{x}),
\end{equation}
where $R$ is the vector of state functions that we introduced in equation~\eqref{eq:stateevolution}. To ensure that the matched density in Section~\ref{subsec:matching} has unit mass, we add an additional macroscopic state function $R_0(\bm{x}) = 1$ to the vector of state functions $\mathbf{R}(\mathbf{x})$ from \eqref{eq:stateevolution} without changing notation. In the context of Monte Carlo simulations, we approximate the continuous probability distribution $\mu$ by a discrete particle ensemble $\mathcal{X} = \left(w_j, \bm{X}_j\right)_{j=1}^J$ with an associated empirical probability distribution $\hat{\mu} = \Sigma_{j=1}^J w_j \delta_{\bm{X}_j}$, and define a discrete version of the restriction operator as
\begin{equation} \label{eq:discreterestrictionoperator}
\mathcal{R}(\mathcal{X}) = \widehat{\mathbb{E}}_{J}[\mathbf{R}(\mathcal{X})] = \frac{1}{J} \sum_{j=1}^J w_j \mathbf{R}(\bm{X}_j).
\end{equation}
The macroscopic state variables at $t^{n} + k\delta t$ are then given by $\mathbf{m}^{n,k}=\mathcal{R}(\mathcal{X}_J^{n,k})$, $k = 0,\dots,K$.

\subsection{Extrapolation} \label{subsec:extrapolation}
To extrapolate the macroscopic state variables over a time interval of size $\Delta t$, we approximate time derivative using a finite difference approximation. Throughout this text, all experiments use linear extrapolation, for which the macroscopic state variables at time $t^{n+1} = t^n+\Delta t$ read
\begin{equation} \label{eq:extrapolation}
\mathbf{m}^{n+1} = \mathbf{m}^n + \frac{\Delta t}{K\delta t}(\mathbf{m}^{n,K}-\mathbf{m}^n).
\end{equation}
Linear extrapolation of the macroscopic state variables mimics a forward Euler step of the unavailable macroscopic model for the macroscopic state variables. Extensions to higher order extrapolation methods are also possible and we refer to \cite{lafitte2017high} for a treatment of higher order projective integration techniques for hyperbolic conservation laws.

\subsection{Matching} \label{subsec:matching} To recover the microscopic level of description from the macroscopic level, we need to find a probability distribution that is consistent with the extrapolated states $\mathbf{m}^{n+1}$. A priori, however, many possible distributions can be consistent with $\mathbf{m}^{n+1}$. To resolve this ill-posedness, we introduce a matching procedure that finds a unique probability distribution that is consistent with the states, while minimizing a dissimilarity functional $\mathcal{I}$ compared to some prior distribution. We choose the prior distribution to be the final distribution $\mu^{n,K}$ at time $t^n+K\delta t$ during the simulation stage in Section~\ref{subsec:mc}. We first introduce the matching operator acting on continuous probability distributions in Section \ref{subsubsec:matching}. In Section \ref{subsubsec:mcmatching} we explain how matching can be implemented on probability distributions represented by weighted particle ensembles.

\subsubsection{The continuous matching operator} \label{subsubsec:matching}
In our approach, matching is formulated as an optimization problem that reads
\begin{equation*}
\mu^{n+1}(\bm{x}) = \underset{\varphi}{\text{arg min}} \ \mathcal{I}(\varphi |\  \mu^{n,K}), \ \ \text{s.t.} \ \ \mathcal{R}(\varphi) = \mathbf{m}^{n+1},
\end{equation*}
with $\mu^{n+1}$ the matched distribution that serves as the initial distribution for the simulation stage in Section \ref{subsec:mc}. Note that the functional $\mathcal{I}$ does not need to be a distance metric. There are many possible matching strategies that come with a different choice for $\mathcal{I}$, see \cite{debrabant2017micro} for more details. In this work, we will use matching in Kullback-Leibler divergence, or relative entropy, which is based on information-theoretic considerations \cite{lelievre2018analysis}. The objective is to minimize
\begin{equation} \label{eq:relativeentropymatching}
\mu^{n+1}(\bm{x}) = \underset{\varphi}{\text{arg min}} \int_G\ln\left(\frac{\varphi(\bm{x})}{\mu^{n,K}(\bm{x})}\right) \varphi(d\bm{x}), \ \ \text{s.t.} \ \ \mathcal{R}(\varphi) = \mathbf{m}^{n+1},
\end{equation}
which has an analytic solution of the form
\begin{equation} \label{eq:matchedanalyticalsolution}
\mu^{n+1}(\bm{x}) = \exp\left(\sum_{k=0}^L \lambda_k R_k(\bm{x}) \right) \mu^{n,K}(\bm{x}),
\end{equation}
where the Lagrange multipliers $\lambda_0, \lambda_1,\dots,\lambda_L$ solve the non-linear system
\begin{equation} \label{eq:nonlinearsystem}
\int_G R_l(\bm{x}) \exp\left(\sum_{k=0}^L \lambda_k R_k(\bm{x})\right) \mu^{n,K}(d\bm{x}) = m_l^{n+1}, \ l = 0,\dots,L.
\end{equation}
In this text, we will employ the Newton-Raphson method to solve the system above to find the Lagrange multipliers \cite{debrabant2017micro}. We give more details in Section~\ref{sec:impl_det}. Minimizing the Kullback-Leibler divergence \eqref{eq:relativeentropymatching} amounts to solving a dual system. 

\subsubsection{Matching with weighted particle ensembles} \label{subsubsec:mcmatching}
In the context of Monte Carlo simulations, we replace the prior distribution $\mu^{n,K}(\bm{x})$ by the particle ensemble $\mathcal{X}^{n,K}$ and the matched distribution by a new particle ensemble $\mathcal{X}^{n+1}$ at time $t^{n+1}=t^n+\Delta t$. We rewrite the matching formula on the ensemble level as
\begin{equation} \label{eq:discretematching}
\mathcal{X}^{n+1} = \mathcal{M}(\mathbf{m}^{n+1},\mathcal{X}^{n,K}),
\end{equation}
where $\mathcal{M}$ is the matching operator that minimizes the Kullback-Leibler divergence \eqref{eq:relativeentropymatching}. 

By the multiplicative nature of matching in Kullback-Leibler divergence \eqref{eq:matchedanalyticalsolution}, sampling from the matched distribution $\mu^{n+1}(x)$ can be performed efficiently by reweighting the particle ensemble $\mathcal{X}^{n,K}$ to represent the matched ensemble $\mathcal{X}^{n+1}$. The new weights read
\begin{equation} \label{eq:particlereweighting}
w_j^{n+1} = \exp\left(\lambda_1R_1(\bm{X}_j^{n,K})+\dots+\lambda_L R_L(\bm{X}_j^{n,K})\right)w_j^{n,K}.
\end{equation}

\subsection{The complete micro-macro acceleration algorithm} \label{subsec:algorithm}
Algorithm \ref{algo:accelerationalgorithm} presents the micro-macro acceleration algorithm in four steps, as defined in Sections~\ref{subsec:mc} to~\ref{subsec:matching}.

\begin{algorithm}
\begin{flushleft}
Assume we have a microscopic ensemble $\mathcal{X}^n = (w_j^n, \bm{X}_j^n)_{j=1}^J$ at time $t^n$, a microscopic step size $\delta t$, the extrapolation step size $\Delta t$ and a number of microscopic steps $K$ such that $K \delta t \leq \Delta t$. Let $T$ be the end time of the simulation and $L$ the number of macroscopic state variables. The algorithm produces the microscopic ensemble $\mathcal{X}^{n+1}$ in four steps:

\vspace{4mm}
(i) \textbf{Monte Carlo simulation}: simulate the microscopic ensemble $\mathcal{X}^n$ over $K$ small inner steps of size $\delta t$
\[
\mathcal{X}^{n, k+1} = S^{n,k}_{\delta t}(\mathcal{X}^{n,k})  \ \ \ \ \ \ \ \ \ \ \ \ \ \  \ \ \ \ \ \  \ \ \eqref{eq:microtimestepper}
\]
		
(ii) \textbf{Restriction}: compute the macroscopic states corresponding to these microscopic ensembles
\[
\mathbf{m}^{n,k} = \mathcal{R}(\mathcal{X}^{n,k}).\ \ \ \ \ \ \ \ \ \ \ \ \ \ \ \ \ \ \ \ \ \ \ \ \ \ \     \eqref{eq:discreterestrictionoperator}
\] 
		
(iii) \textbf{Extrapolation}: approximate the macroscopic state variables over a larger time interval $\Delta t \geq K\delta t$ with linear extrapolation
\[
\mathbf{m}^{n+1} = \mathbf{m^n} + \frac{\Delta t}{K \delta t}\left(\mathbf{m}^{n,K} - \mathbf{m}^n \right) \ \ \ \ \eqref{eq:extrapolation}
\]
		
(iv) \textbf{Matching}: compute the new particle ensemble consistent with $\mathbf{m}^{n+1}$ by reweighting the prior ensemble $\mathcal{X}^{n,K}$ using the matching operator
\[
\mathcal{X}^{n+1} = \mathcal{M}(\mathbf{m}^{n+1},\mathcal{X}^{n,K}) \ \ \ \ \ \ \ \ \ \ \ \ \ \ \ \ \ \eqref{eq:discretematching}
\]
and advance time with $\Delta t$ until the end time $T$ is reached.
\end{flushleft}
\caption{Micro-macro acceleration}
\label{algo:accelerationalgorithm}
\end{algorithm}

\section{Practical implementation details} \label{sec:impl_det}
To implement the micro-macro acceleration algorithm in practice, several issues need to be addressed. In this Section, we focus on three aspects. First, in Section~\ref{subsec:newtonraphson}, we formulate an efficient implementation of the Newton-Raphson scheme to solve the non-linear system~\eqref{eq:nonlinearsystem} using the weighted particle ensemble available at the end of the simulation stage. 

Second, due to the particle reweighting~\eqref{eq:particlereweighting}, some particles can end up with very large or very small weights. Such a variation in the weights can render the approximation of the macroscopic state variables in~\eqref{eq:discreterestrictionoperator} inaccurate as some particles with small weights are essentially ignored. In Section~\ref{subsec:resampling}, we go deeper into the problem of detecting a large variation in the weights. To alleviate this variation, we use a resampling strategy that randomly duplicates some particles and removes others in such way that, after resampling, all resulting particles have equal weights~\cite{hol2006resampling}.

Finally, in practice it can happen that the extrapolated macroscopic states $\textbf{m}^{n+1}$ fall outside the domain of the matching operator $\mathcal{M}$, i.e., there is no probability distribution that is consistent with $\textbf{m}^{n+1}$. We call such a phenomenon a \textit{matching failure}. An efficient way of dealing with matching failures is decreasing the extrapolation step size such that the extrapolated states deviate less from the already available states in the simulation stage. The detection of matching failures and the adaptive time stepping strategy are discussed in Section~\ref{subsec:match_fails}.

\subsection{An efficient formulation of the Newton-Raphson scheme} \label{subsec:newtonraphson}
To compute the matched distribution $\mu^{n+1}$, we need to find the Lagrange multipliers $\bm{\lambda}$ that solve the dual system of non-linear equations $g(\bm{\lambda}) = 0$ \eqref{eq:nonlinearsystem}, where we define
\begin{equation}
    g(\bm{\lambda})_l = m_l - \int_G R_l(\mathbf{x}) \exp \left(\sum_{m=0}^L \lambda_m R_m(\mathbf{x}) \right) \mu^{n,K}(\mathbf{x}) \ d\mathbf{x},\ \ l = 0,\dots,L.
\end{equation}
Suppose we have the Lagrange multipliers at the $k$-th iteration of the Newton-Raphson scheme, $\bm{\lambda}^{(k)} = \left( \lambda_0^{(k)}, \dots, \lambda_L^{(k)}\right)^T$. One Newton-Raphson iteration for the next set of Lagrange multipliers $\bm{\lambda}^{(k+1)}$ for the dual problem \eqref{eq:nonlinearsystem} reads
\begin{equation} \label{eq:newtonraphson}
\bm{\lambda}^{(k+1)} = \bm{\lambda}^{(k)} - (\nabla g(\bm{\lambda}^{(k)}))^{-1}  g(\bm{\lambda}^{(k)}),
\end{equation}
where the gradient of $g(\bm{\lambda})$ is given by
\begin{equation*}
\nabla g(\bm{\lambda})_{k,l} = - \int_G R_k(\mathbf{x}) R_l(\mathbf{x}) \exp\left(\sum_{m=0}^L \lambda_m R_m(\mathbf{x})\right)  \mu^{n,K}(\mathbf{x}) \ d\mathbf{x}, \ k,l=0,\dots,L
\end{equation*}
In our setting, the prior probability distribution $\mu^{n,K}$ is only available as an empirical measure $\hat{\mu}^{n,K}$. We can approximate the function $g$ and its Jacobian using the prior particle ensemble $\mathcal{X}^{n,K} = \left(w_j^n, \bm{X}_j^{n,K}\right)_{j=1}^J$ as
\begin{align}
\begin{split}
\hat{g}(\bm{\lambda})_l &= m_l - \sum_{j=1}^J R_l\left(\mathbf{X}^{n,K}_j\right) \exp \left(\sum_{m=0}^L \lambda_m R_m\left(\mathbf{X}^{n,K}_j\right) \right) w_j^{n} ,\ \ l = 0,\dots,L, \\
\nabla \hat{g}(\bm{\lambda})_{k,l} &= - \sum_{j=1}^J R_k\left(\mathbf{X}^{n,K}_j\right) R_l\left(\mathbf{X}^{n,K}_j\right) \exp\left(\sum_{m=0}^L \lambda_m R_m\left(\mathbf{X}_j^{n,K}\right)\right) w_j^{n}, \ \ k,l=0,\dots,L.
\end{split}
\end{align}
Note that since the prior particle ensemble $\mathcal{X}^{n,K}$ is fixed, the solution to the system $\hat{g}(\bm{\lambda}) = 0$ is also deterministic.

The cost of one Newton-Raphson iteration increases linearly with the number of particles $J$, making matching the most expensive part of Algorithm \ref{algo:accelerationalgorithm}. There exist other numerical techniques, such as a sparse approximation of the Jacobian \cite{lucia1983explicit} or the BFGS algorithm \cite{martinez2000practical} that could speedup matching. These methods usually require more iterations to converge to the correct multipliers, but the iterations require fewer computations. We use the Newton-Raphson procedure in all the numerical experiments in this manuscript.

As initial value for \eqref{eq:newtonraphson}, we take $\bm{\lambda}^{(0)} = 0 \in \mathbb{R}^{L+1}$. When there is no extrapolation, i.e. $\Delta t = \delta t$, the extrapolated macroscopic state variables are consistent with the prior distribution $\mu^{n,K}$, i.e., $\mathcal{R}(\mu^{n,K}) = \mathbf{m}^{n+1}$. In this case, the optimal Lagrange multipliers equal 0 so that equations~\eqref{eq:newtonraphson} and \eqref{eq:nonlinearsystem} are already solved by the initial guess. As stopping criterion we require that extrapolated moments are met, within a tolerance
\[
\norm{\hat{g}(\bm{\lambda}^k)}_2 < \text{TOL}
\]
where the tolerance is usually $\text{TOL} = 10^{-9}$. In practice, we also impose a maximum number of Newton-Raphson iterations due to the possibility of \textit{matching failures}, which we will explain further in Section~\ref{subsec:match_fails}.

\subsection{Particle resampling}~\label{subsec:resampling}
Besides the problem of matching failures, there is also the issue of increasing variance between particle weights. By the multiplicative form of the matching procedure \eqref{eq:matchedanalyticalsolution} and particle reweighting \eqref{eq:particlereweighting}, some Monte Carlo replicas may end up with a high weight after a few steps, while some may have a low weight. Particles with a small weight will be neglected in the discrete approximation of the macroscopic state
variables~\eqref{eq:discreterestrictionoperator}, while a few particles with large weights will dominate the approximation. Hence, the discrete approximation to the macroscopic state variables becomes prone to statistical errors. This is in fact a well-known issue in many particle filter methods \cite{arulampalam2002tutorial}.

To resolve the large variation in particle weights, we resample the particles from time to time. The technique we use in this paper is \textit{stratified resampling} \cite{hol2006resampling,debrabant2017micro}. The idea is to randomly discard and duplicate each of the particles in such way that the distribution of the particles is unchanged, and all particles have an equal weight $1/J$ after resampling. In practice, one divides the unit interval $[0,1]$ in $J$ equal pieces and generates a uniformly distributed random number in each sub-interval
\[
u_k = \frac{(k-1)+\tilde{u}_k}{J}, \ \ \tilde{u}_k \sim \mathcal{U}(0,1).
\]
One then duplicates each particle $\bm{X}_j^{n+1}$, $n_j$ times where $n_j$ is chosen by counting how many bins of equal size $1/J$ in $[0,1]$ the weight $w_j^{n+1}$ covers \cite{douc2005comparison}. Mathematically, we thus take $n_j = \# \{u_k | u_k \in [ \sum_{i=1}^{j-1} w_i^{n+1}, \sum_{i=1}^j w_i^{n+1}] \}$. If $n_j=0$ we discard particle $\bm{X}_j^{n+1}$ and if $n_j >1$ we duplicate the particle into more particles. When $n_j=1$ then the particle remains and we simply adjust its weight to $1/J$. The above procedure ensures that the particle distribution remains unchanged, in expectation \cite{douc2005comparison}.

To determine whether the variation between the weights is too large, we  propose to compute the discrete relative entropy of the weights compared to uniform weights $1/J$,
\[
\sum_{j=1}^J w_j^{n+1} \ln(J w_j^{n+1}) \in [0, \ln J],
\]
which is zero when all weights are equal. We propose, as a heuristic, to resample when this relative entropy exceeds a fixed value of $\ln(J)/10$. We choose this value to ensure that the variation in the weights is not too high and prevent inaccurate computations of the macroscopic states. In practice, we compute the relative entropy of the weights every five time steps of the micro-macro acceleration algorithm, and we resample the weights if the entropy exceeds the chosen maximal value.

\subsection{An adaptive time-stepping strategy for matching failures} \label{subsec:match_fails}
As we mentioned in the introduction of this Section, a problem arises when the extrapolated macroscopic state variables $\mathbf{m}^{n+1}$ fall outside the domain of the matching operator $\mathcal{M}$. In this case, there is no probability distribution that is consistent with the extrapolated states, and the micro-macro acceleration algorithm fails. We call such behavior a \textit{matching failure} \cite{lelievre2018analysis,debrabant2017micro}. When the extrapolation step $\Delta t$ is smaller, the extrapolated macroscopic state variables will lie closer to the macroscopic state variables of the prior distribution, for which there exists a probability distribution consistent with the prior state variables, namely the prior distribution itself. We can thus expect that, the smaller the extrapolation step, the more likely we will find a probability distribution consistent with $\mathbf{m}^{n+1}$. In practice, we thus perform micro-macro acceleration with \textit{adaptive time stepping} based on matching failures. If we detect a matching failure, we drastically decrease $\Delta t$ by half, and if matching succeeds we cautiously increase $\Delta t$ again by a factor $1.2$.

The question remains how we can detect a matching failure in practice. When there is no probability distribution consistent with the extrapolated state variables $\mathbf{m}^{n+1}$, the Newton-Raphson solver will also fail to find Lagrange multipliers that solve the non-linear system~\eqref{eq:nonlinearsystem}. We thus only need to monitor the convergence of the Newton-Raphson solver to determine whether there is a matching failure or not. 
We impose a maximum number of iterations for the non-linear solver. If this maximum number is reached before convergence, we consider matching to have failed.

Imposing a maximum number of iterations before lowering the extrapolation time step size also has another beneficial effect. For larger $\Delta t$, the extrapolated state variables will lie far from the prior macroscopic state variables so that the Newton-Raphson method will need more iterations to convergence compared to lower extrapolation step sizes. Taking two time steps with a lower $\Delta t$ could speed up the algorithm compared to one step with twice the extrapolation time step. In the light of performance, we can hence take the maximum number of Newton-Raphson iterations to be small. In this text, we allow 6 Newton-Raphson iterations before lowering the time step.

\section{The influence of the choice of state variables} \label{sec:states}
In this Section, we present two examples from molecular dynamics: FENE-dumbbells (Section~\ref{subsec:fene}) and a three-atom molecule (Section~\ref{subsec:threeatom}). In both examples, we investigate the accuracy of micro-macro simulation for different choices of macroscopic state variables. In Section~\ref{subsec:fene}, we introduce the model of FENE-dumbbells and look at three different hierarchies of macroscopic state variables to approximate the quantities of interest. Section~\ref{subsec:threeatom} introduces and discusses the three-atom molecule, where we compare the accuracy of micro-macro acceleration to existing approximate macroscopic models for slow quantities of interest. We will show that even with some macroscopic state variables that are not entirely slow, the micro-macro acceleration method can generate more accurate approximations to the exact dynamics than these approximate models. Finally, we present a short overarching discussion on the effect of the macroscopic state variables in Section~\ref{subsec:conlcusomoments}.

\subsection{Three state hierarchies for FENE-dumbbells} \label{subsec:fene}
FENE stands for `Finitely Extensible Non-linear Elastic', and FENE-dumbbells model represents a dilute polymer solution in a solvent\cite{laso1993calculation}. The state variable $X(t) \in \mathbb{R}$ denotes the end-to-end vector between two beads of the polymer chain, which are connected by a spring. As the polymers move through the solvent, they experience Stokes drag due to the velocity gradient $\kappa(t)$ of the ambient solvent, Brownian motion $W(t)$ due to collisions with solvent molecules and a spring force
\begin{equation}
F: B(b) \to \mathbb{R}: x \mapsto \frac{b^2x}{b^2-x^2},
\end{equation}
due to intramolecular interactions. The set $B(b)$ is the open ball of radius $b$ around the origin. The diffusion process for FENE-dumbbells reads
\begin{equation} \label{eq:feneprocess}
dX_t = \left(\mathbf{\kappa}(t)X_t - \frac{1}{2\text{We}}F(X_t)\right)dt + \frac{1}{\sqrt{\text{We}}}dW_t.
\end{equation}
where $\text{We}> 0$ is the Weissenberg number, the ratio of the characteristic relaxation time of the polymers in the solvent to the characteristic time of the solvent. We refer to \cite{bird1987dynamics} for a derivation of SDE~\eqref{eq:feneprocess} and to \cite{jourdain2003mathematical,jourdain2004existence} for the existence of a global solution to the SDE.

\begin{figure}
\centering
\includegraphics[width=0.45\textwidth]{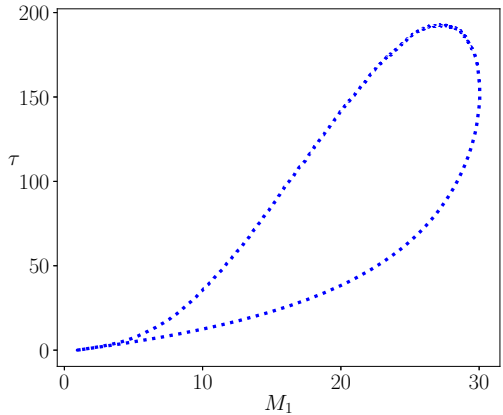}
\caption{The $\tau-M_1$ phase diagram for the FENE process \eqref{eq:feneprocess} with the time-dependent velocity field, computed with the Euler-Maruyama method with time step $\delta t = 2 \cdot 10^{-4}$ up to 4 seconds.}
\label{fig:fene_hysteresis}
\end{figure}

In practice, the FENE-dumbbells process is usually coupled to the macroscopic Navier-Stokes equations that models the evolution of the solvent. This coupling happens through the Stokes drag $\kappa(t)$ and by adding a non-Newtonian term in the form of the stress tensor \cite{masmoudi2008well}
\begin{equation*}
\mathbf{\tau}(t) = \frac{1}{\text{We}} \left( \mathbb{E}[X_t F(X_t)] - 1 \right).
\end{equation*}
In the following experiments, we will use a 1-dimensional FENE-dumbbells process, with parameters $W_{\text{e}}=1, \ b = 7$. We also define the ambient velocity gradient as in \cite{keunings1997peterlin}
\[
\kappa(t) = 100t(1-t)e^{-4t}.
\]
For this choice of the velocity gradient, one can show that hysteresis occurs between the stress tensor and the first even moment, $M_1 = \mathbb{E}[X(t)^2]$ as shown on the $\tau-M_1$ phase diagram in Figure~\ref{fig:fene_hysteresis}. The hysteric effect is due to the interplay between the non-linearities in the FENE model \eqref{eq:feneprocess} and the dispersity between the elongations of the individual polymers. We refer to \cite{sizaire1999hysteretic} for more details on this hysteric effect. 

The objective of the numerical experiments in this Section is to study the accuracy of micro-macro acceleration in approximating this hysteretic curve, and the effect of the choice of the macroscopic state variables on this accuracy. We propose three different hierarchies of macroscopic state variables and study the accuracy of the resulting approximations by micro-macro acceleration.
The micro-macro acceleration algorithm has already been successfully applied to the FENE-dumbbells diffusion process \cite{debrabant2017micro}, but a study of the accuracy of different hierarchies of macroscopic state variables has not yet been performed. The authors of \cite{samaey2011numerical} have, however, already carried out a similar study in the equation-free context.

In the following experiments, we will use the same three hierarchies of macroscopic state variables as \cite{samaey2011numerical}, and investigate the accuracy of micro-macro acceleration on the hysteretic curve and the stress tensor, as a function of the hierarchy and the number of macroscopic state variables per hierarchy. The families of macroscopic state variables are:

\begin{itemize}
\item Hierarchy 1: Use the first $L$ even moments of the diffusion process. By symmetry of the spring force, the odd moments of the FENE process vanish, so there is no need to use these states for matching. The macroscopic state variables read
\begin{equation*}
\bm{m}(t) = \left( \mathbb{E}[X_t^2], \dots,\mathbb{E}[X_t^{2L}] \right), \ L > 0.
\end{equation*}
\item Hierarchy 2: Replace the final state variable $\mathbb{E}[X^{2L}]$ from strategy 1 with the stress tensor $\tau$ itself. Since we want to approximate the evolution stress tensor, adding the latter as a state variable could improve the accuracy of micro-macro acceleration. The second hierarchy of state variables hence reads
\begin{equation*}
    \bm{m}(t) = \left( \mathbb{E}[X_t^2], \dots,\mathbb{E}[X_t^{2(L-1)}], \tau(t) \right), \ L > 0.
\end{equation*}
\item Hierarchy 3: Start with $M_1$ and add terms that pop up in the Taylor expansion of $\tau$ using It{\^o}'s lemma. For each of those new terms that pop up in the evolution equation for $M_1$, write down their evolution equation as well and keep adding new terms as macroscopic state variables. This way, the first four state variables that pop up are
\begin{equation} \label{eq:strategy3}
\bm{m}(t) = \left(\mathbb{E}[X_t^2],  \mathbb{E}\left[\frac{X_t^2}{1-X_t^2/b^2}-1\right],  \mathbb{E}\left[\frac{X_t^2}{(1-X_t^2/b^2)^2}\right], \mathbb{E}\left[\frac{X_t^2}{(1-X_t^2/b^2)^3}\right] \right).
\end{equation}
We refer to Appendix~\ref{app:fene} for the derivation of the above macroscopic state variables.
\end{itemize}

\begin{figure}
\centering
\includegraphics[width=0.8\textwidth]{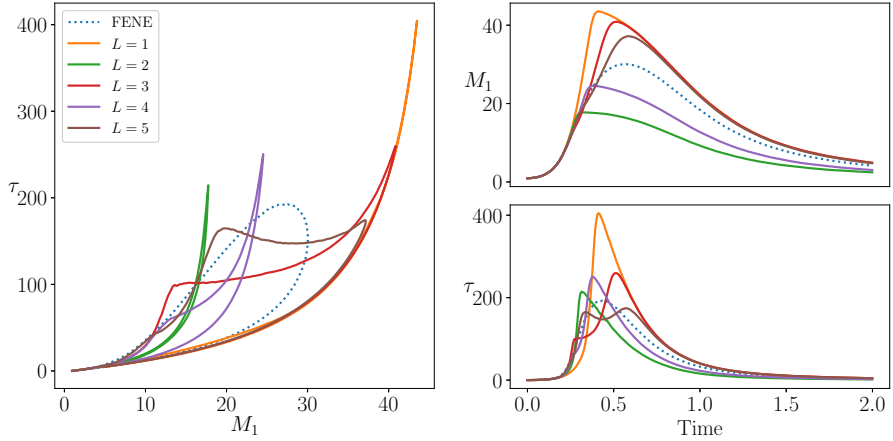}
\caption{The $\tau$-$M_1$ phase diagram for the first hierarchy of macroscopic state variables (left) and the evolution of $M_1$ (top right) and of the stress tensor $\tau$ (Bottom right). Only for $L=5$, the first hierarchy starts to approximate the hysteresis curve somewhat, although there still is room for improvement.}
\label{fig:fene_strategy1}
\end{figure}

\begin{figure}
	\centering
	\includegraphics[width=0.8\textwidth]{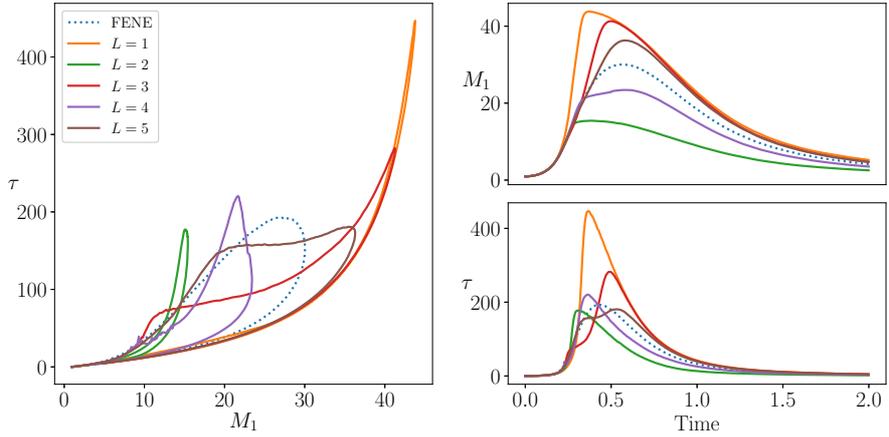}
	\caption{The $\tau$-$M_1$ phase diagram for the second hierarchy of macroscopic state variables (left) and the evolution of $M_1$ (top right) and of the stress tensor $\tau$ (Bottom right). The hysteresis curve is somewhat better approximated, but the result is still not convincing.}
	\label{fig:fene_strategy2}
\end{figure}

\begin{figure}
\centering
\includegraphics[width=0.8\textwidth]{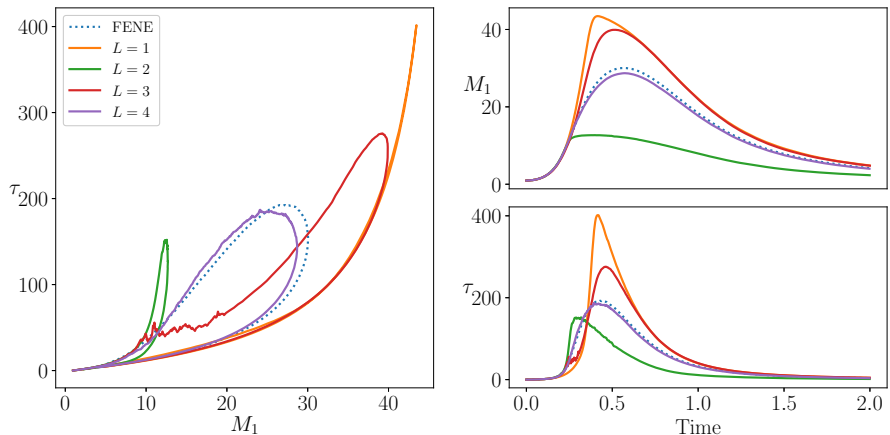}
\caption{The $\tau$-$M_1$ phase diagram for the hierarchy hierarchy of macroscopic state variables (left) and the evolution of $M_1$ (top right) and of the stress tensor $\tau$ (Bottom right). The third hierarchy of states already gives much closer approximations to the hysteresis curve when $L=4$, compared to the two previous hierarchies that require at least $L=5$ states for a decent approximation.}
\label{fig:fene_strategy3}
\end{figure}

For the numerical experiments, we compute the evolution of the stress tensor and the first even moment $M_1$ and also plot the hysteric curve in the $\tau-M_1$ phase space. We will use up to $L=5$ macroscopic state variables for the first two hierarchies and up to $L=4$ state variables for the third hierarchy~\eqref{eq:strategy3}. We take a maximal time step $\Delta t = 5 \delta t$ for the extrapolation time step with adaptive time stepping, since the expected gain is quite small due to a small relative time-scale separation in the FENE model. Furthermore, we choose $\delta t = 2 \cdot 10^{-4}$ for all experiments and we use $N=5 \cdot 10^4$ particles during the Monte Carlo simulations. The numerical results are summarized in Figures \ref{fig:fene_strategy1} to \ref{fig:fene_strategy3}, for each of the three hierarchies of state variables, respectively. 

With the first hierarchy of macroscopic state variables in Figure \ref{fig:fene_strategy1}, micro-macro acceleration does not capture the hysteresis curve well for $L=1,\dots, 4$. For $L=5$ in brown, the approximation to the exact hysteresis curve is somewhat better, but there is still a lot of room for improvement.

The second hierarchy of state variables, which includes the stress tensor, improves the approximation of micro-macro acceleration in Figure \ref{fig:fene_strategy2}. The green and purple lines enclose a larger area and the brown line approximates the hysteresis curve slightly better. There is however still no reasonable fit of the exact hysteresis curve when $L=5$.

The third hierarchy of macroscopic states yields the best results. For $L=4$, micro-macro acceleration (in purple) follows the exact hysteresis curve closely. With the third hierarchy of state variables, we hence perform less computational work while obtaining a better approximation to the exact stress tensor and hysteresis curve. A priori selecting the set of macroscopic state variables has a profound impact on the accuracy and efficiency of micro-macro acceleration. This conclusion is very similar to the one in the equation-free context~\cite{samaey2011numerical}.

\subsection{Reaction coordinates of a three-atom molecule} \label{subsec:threeatom}
\begin{figure}
	\centering
	\includegraphics[width=0.25\textwidth]{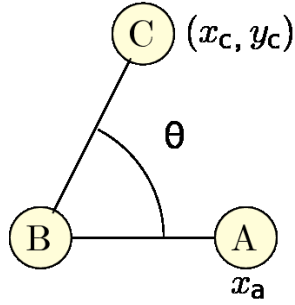}
	\caption{The three-atom molecule. Atom A is restricted to the $x$-axis and B is fixed at the origin.}
	\label{fig:tri_atom_molecule}
\end{figure}

The three-atom molecule is a simple three-dimensional molecular system, where the bonds between the individual atoms vibrate at a high frequency, compared to the frequency of global conformational changes of the molecule. Figure \ref{fig:tri_atom_molecule} depicts the simple molecule. To remove some degrees of freedom, atom A is restricted to move on the $x$-axis with coordinates $(x_a, 0)$ and atom B is fixed at the origin. We define the coordinates for atom C as $(x_c, y_c)$. The individual atoms are influenced by a potential that describes the interactions between the atoms, and also by collisions with the ambient solvent, modeled via Brownian motion. The magnitude of these collisions is proportional to the inverse temperature of the solvent.

The evolution of the atoms is given by an overdamped Langevin dynamics of the form\cite{legoll2010effective}
\begin{equation} \label{eq:triatommolecule}
\begin{cases}
dx_a = -\frac{\partial V}{\partial x_a} dt + \sqrt{2 \beta^{-1}}dW_{x_a} \\
dx_c = -\frac{\partial V}{\partial x_c} dt + \sqrt{2 \beta^{-1}} dW_{x_c} \\
dy_c = -\frac{\partial V}{\partial y_c} dt + \sqrt{2 \beta^{-1}} dW_{y_c},
\end{cases}
\end{equation}
where $\beta$ is the inverse temperature. The potential energy that governs the drift of the individual atoms is given as
\begin{equation} \label{eq:triatompotential}
V(x_a, x_c, y_c) = \frac{1}{2\varepsilon}(x_a-l_{\text{eq}})^2 + \frac{1}{2\varepsilon}(\sqrt{x_c^2+y_c^2}-l_{\text{eq}})^2 + \frac{k}{2}((\theta - \theta_{\text{saddle}})^2 - \delta \theta^2),
\end{equation}
where $l_{\text{eq}}$ is the equilibrium length for the bonds and $\theta$ is the angle between atoms A and C. As we can see from the expression of the potential energy, the atom has two stable conformations where $\theta = \theta_{\text{saddle}} \pm \delta \theta$. Due to the time-scale separation in the potential energy, expressed by the small scale parameter $\varepsilon \ll 1$, the bonds between atoms A and B and atoms B and C vibrate quickly relative to the bimodal behavior of the angle $\theta$, which is the variable of interest in this system. For the experiments in this section, we will use the same parameter values as in \cite{legoll2010effective}, i.e., $\beta = 1, \ \varepsilon = 10^{-3}, \ k = 208, \ \theta_{\text{saddle}} = \pi/2$ and $\delta \theta = \theta_{\text{saddle}} - 1.1187$.

In Section~\ref{subsubsec:eff} we describe how an approximate model, or an effective dynamics model can be constructed for a quantity of interest, or reaction coordinate, when the time-scale separation is large. In Section~\ref{subsubsec:choicer} we choose two reaction coordinates and discuss the accuracy of the effective dynamics. In Section~\ref{subsubsec:mMtri}, we apply micro-macro acceleration to the three-atom molecule with macroscopic state variables closely related to the quantity of interest (the angle $\theta$) and look at its efficiency gain. Finally, in Section~\ref{subsubsec:badmac}, we look at the accuracy of micro-macro acceleration with macroscopic state variables that are loosely coupled to $\theta$.

\subsubsection{Effective dynamics} \label{subsubsec:eff}
When the parameter $\varepsilon$ is small, a direct simulation of \eqref{eq:triatommolecule} becomes prohibitively costly. In many situations in molecular dynamics, however, we are usually not interest in the evolution of the complete stochastic system, but rather in the evolution of a low-dimensional set of reaction coordinates. For the three-atom molecule~\eqref{eq:triatommolecule}, we consider a scalar reaction coordinate $\xi = \xi(x_a,x_c,y_c)$. In \cite{legoll2010effective}, the authors propose to use a so-called effective dynamics model, or an approximate macroscopic model for the reaction coordinate $\xi$. The effective dynamics is based on the conditional expectation of the state variable $X=(x_a, x_c, y_c)$ given a value of the reaction coordinate $\xi$, and reads
\begin{equation} \label{eq:effecitvedynamics}
d\xi = b(\xi)dt + \sqrt{2 \beta^{-1}} \sigma(\xi)dW_{\xi}.
\end{equation}
The drift $b(\xi)$ and diffusion term $\sigma(\xi)$ are defined by the conditional expectations
\begin{align} \label{eq:effectivedynamicsfunctions}
b(z) &= \mathbb{E}_{\Psi_{\infty}}[-\nabla V(X) \cdot \nabla \xi(X) + \beta^{-1}\Delta \xi(X) \ | \ \xi(X) = z] \\
\sigma(z)^2 &= \mathbb{E}_{\Psi_{\infty}}[\ \norm{\nabla \xi(X)}^2 \ | \ \xi(X) = z],
\end{align}
where $\Psi_{\infty} = \frac{1}{Z}\exp\left(\beta^{-1}V\right)$ is the invariant measure of \eqref{eq:triatommolecule} with $Z$ the normalization constant. We refer to \cite{legoll2010effective} for convergence properties of the effective dynamics \eqref{eq:effecitvedynamics}. 

Analytic expressions for the functions $b(z)$ and $\sigma(z)$ are usually hard or impossible to derive. Therefore, we pre-compute the functions on a grid of $z$-values and interpolate between these points when an intermediate value is required \cite{legoll2010effective}.

\subsubsection{Choice of reaction coordinates} \label{subsubsec:choicer}
\begin{figure}
\centering
\begin{subfigure}[b]{0.5\textwidth}
\centering
\includegraphics[width=0.8\textwidth]{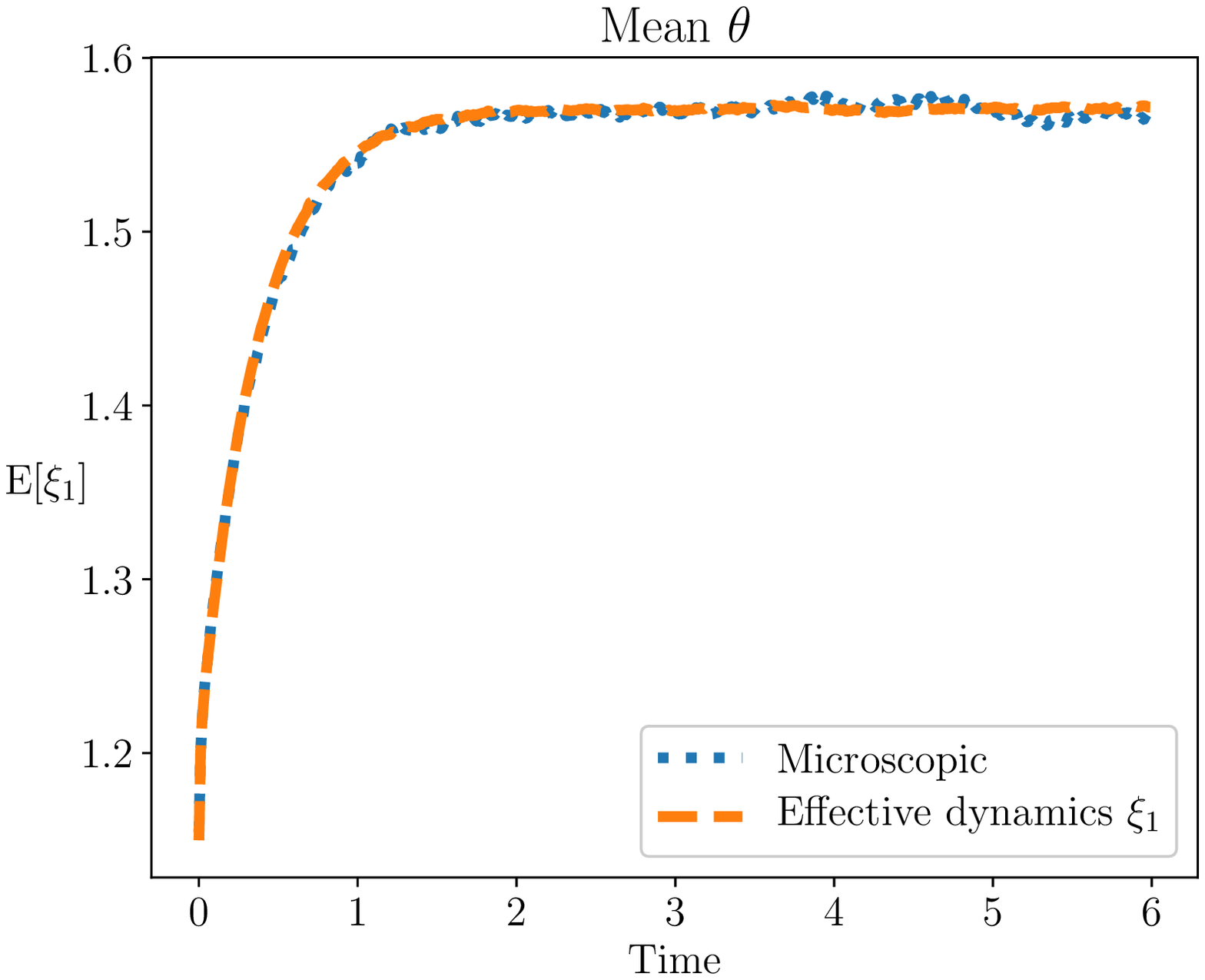}
\end{subfigure}%
\begin{subfigure}[b]{0.5\textwidth}
\centering
\includegraphics[width=0.8\textwidth]{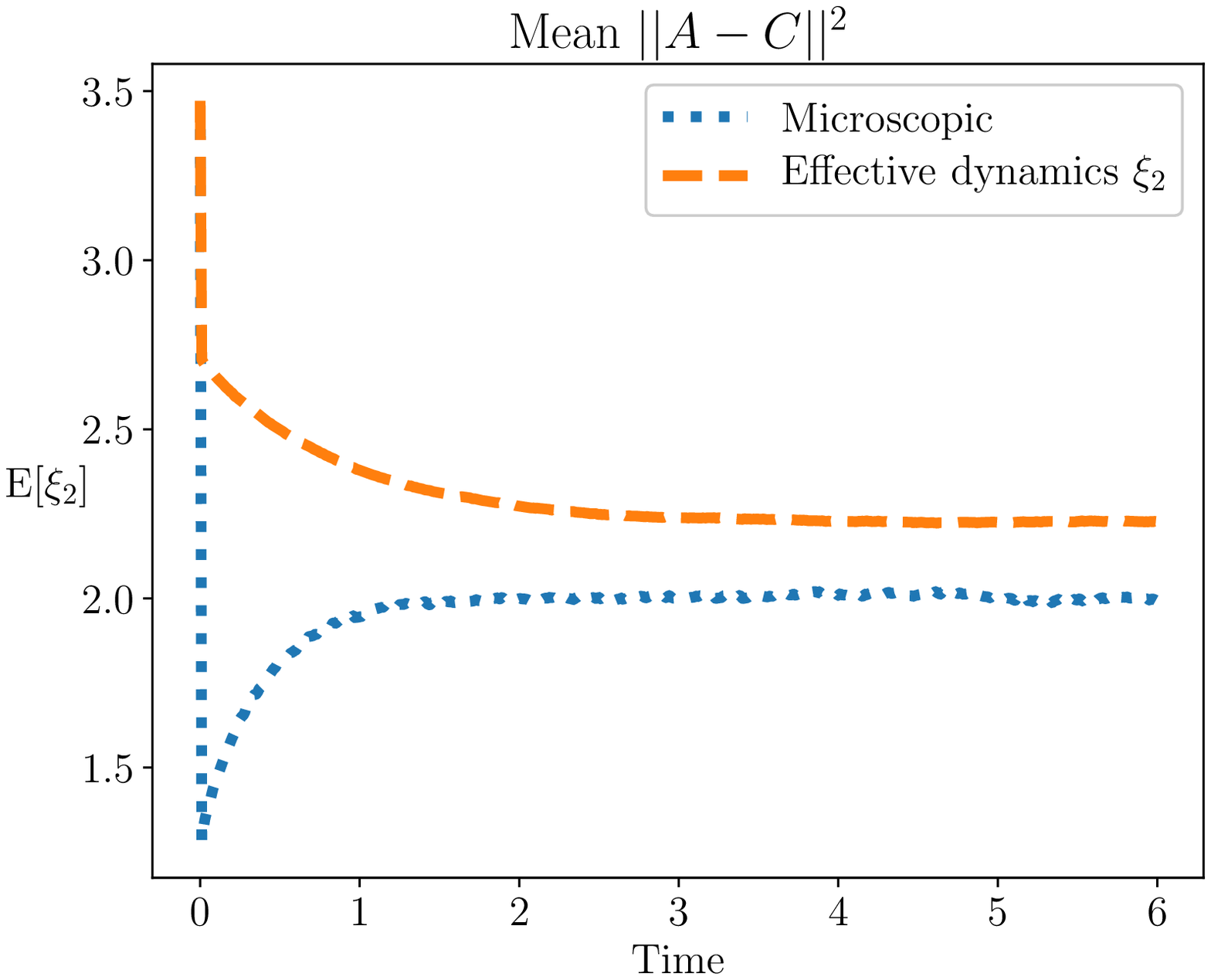}
\end{subfigure}
\caption{Left: The exact evolution of the mean of $\theta$ given by \eqref{eq:triatommolecule}, computed with the Euler-Maruyama method in blue and the effective dynamics of reaction coordinate $\xi_1$ in orange. Right: a similar plot of the evolution of the mean of $\norm{A-C}^2$ with the effective dynamics based on $\xi_2$ in orange. Reaction coordinate $\xi_2$ makes a modeling error while $\xi_1$ does not.}
\label{fig:triatomapproximatemodels}
\end{figure}

As we mentioned in the introduction of this Section, we are mostly interested in the possible conformations that the molecule can reside in. To this end, we define a first reaction coordinate $\xi_1$ to be the angle $\theta$. Alternatively, we define a reaction coordinate $\xi_2$ as the $L_2$-distance squared between the two outer atoms A and C, i.e.,
\begin{align*}
\xi_1 &= \theta, \\
\xi_2 &= \norm{A-C}^2.
\end{align*}
The effective dynamics evolution equation for $\theta$ is readily given by the potential energy function \eqref{eq:triatommolecule}, while for the second reaction coordinate, we compute the drift and diffusion terms, $b(z)$ and $\sigma(z)$ on uniform grid of $z$-values between 0 and 5.
The reaction coordinates seem very similar, as they both can capture the bimodal nature of the dynamics of $\theta$. However, the accuracy of the resulting effective dynamics differ considerably. The numerical results for $\xi_1$ and $\xi_2$, as reproduced from \cite{legoll2010effective}, are summarized in Figure \ref{fig:triatomapproximatemodels}. The microscopic time step is $\delta t = 10^{-3}$ for all simulations and we use $N=5 \cdot 10^4$ Monte Carlo replicas.

The effective dynamics based on $\xi_1$ follows the exact dynamics of the angle $\theta$ very well. On the other hand, the effective dynamics with $\xi_2$ does not follow the exact dynamics of $\norm{A-C}^2$ at all. There is a large steady-state error between the effective dynamics and the exact microscopic dynamics. The authors of \cite{legoll2010effective} attribute this error to the fact that $\xi_2$ is not orthogonal to the fast dynamics, while $\xi_1$ is. The a priori choice of which reaction coordinate to use, is thus of great importance.

\subsubsection{Efficiency gain by micro-macro acceleration} \label{subsubsec:mMtri}
\begin{figure}
\centering
\begin{subfigure}[b]{0.5\textwidth}
\centering
\includegraphics[width=0.85\textwidth]{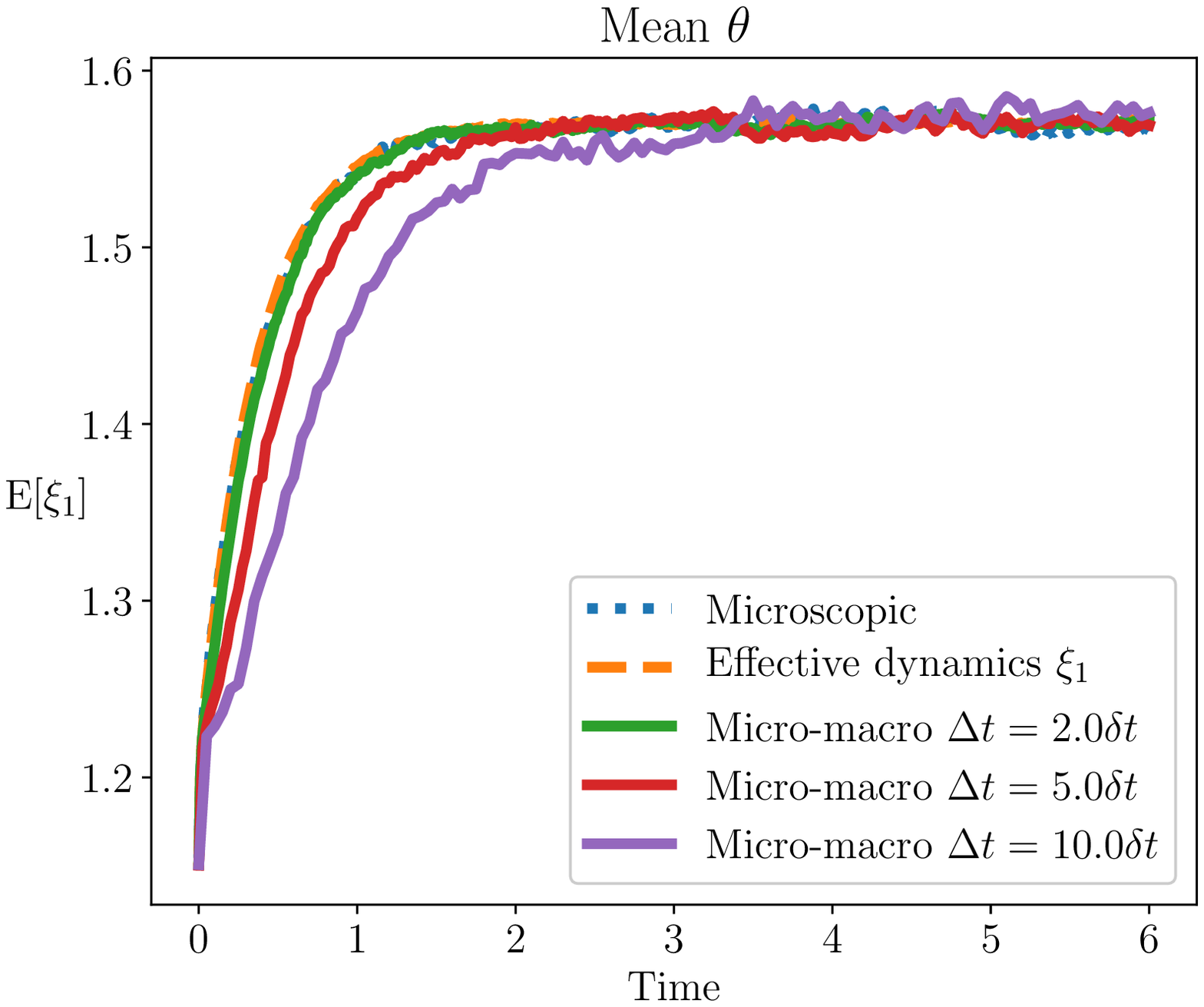}
\end{subfigure}%
\begin{subfigure}[b]{0.5\textwidth}
\centering
\includegraphics[width=0.85\textwidth]{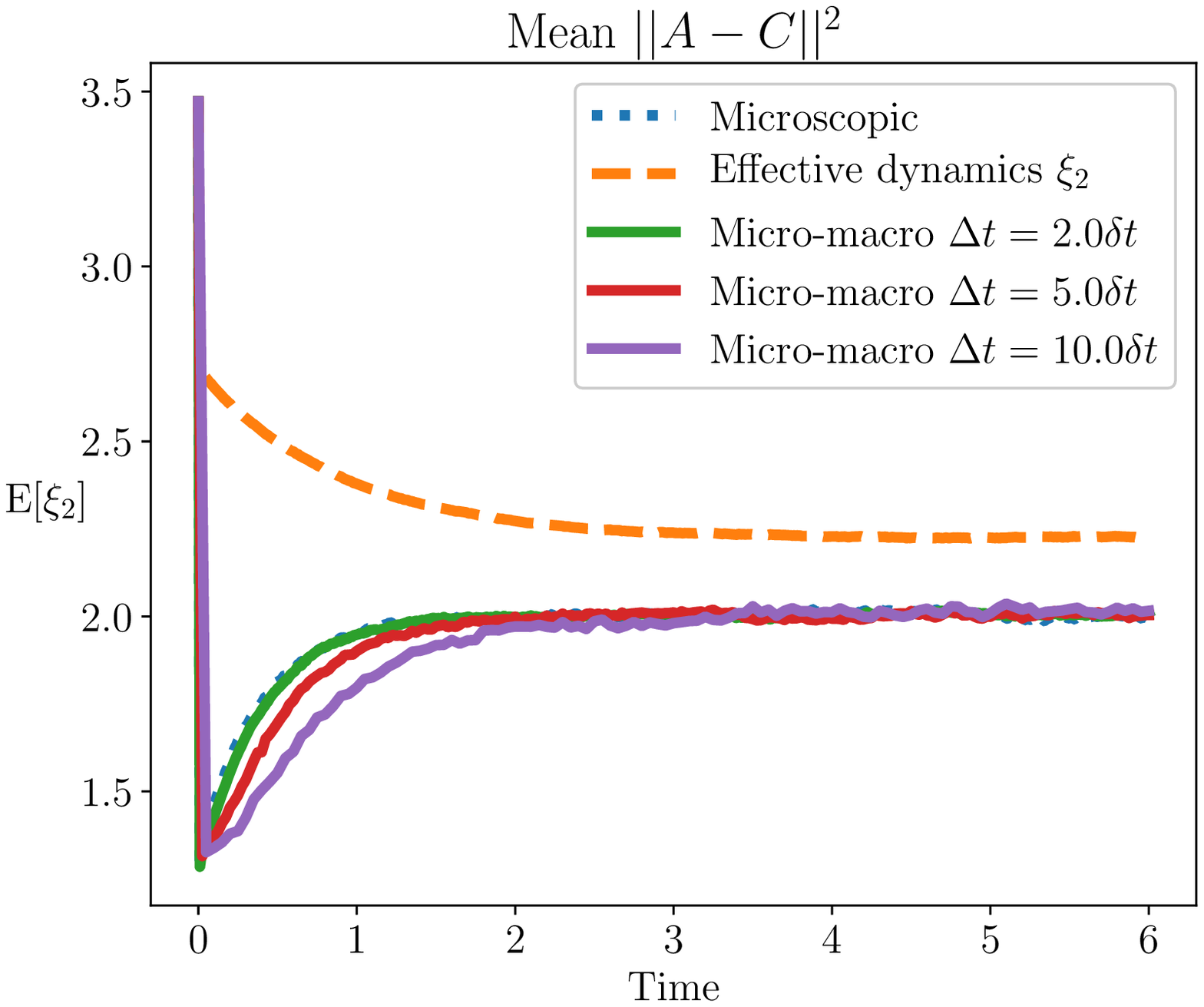}
\end{subfigure}
\caption{Left: same plot as Figure \ref{fig:triatomapproximatemodels} (left). The solid lines represent the solution obtained by micro-macro acceleration with only extrapolating the mean of $\xi_1$, for different extrapolation step sizes. Right: same plot as  Figure \ref{fig:triatomapproximatemodels} (right) where we extrapolate the mean of $\xi_2$ as a macroscopic state. Micro-macro acceleration removes the modeling error of the effective dynamics of $\xi_2$.}
\label{fig:triatommicro-macro}
\end{figure}

The choice of a good set of reaction coordinates to use can be very hard in more complex applications \cite{legoll2010effective}. Micro-macro acceleration can reduce this difficulty. We now perform the same numerical experiments as above with the micro-macro acceleration method, and study when micro-macro acceleration can improve on the results of the effective dynamics. In the following experiments, we extrapolate the mean of either $\xi_1$ or $\xi_2$ as a macroscopic state variable, for different extrapolation step sizes. The simulation results are summarized in Figure \ref{fig:triatommicro-macro}, with reaction coordinate $\xi_1$ on the left and $\xi_2$ on the right.

On the left plot of Figure \ref{fig:triatommicro-macro}, we see that micro-macro acceleration closely follows the exact dynamics of the first reaction coordinate, for different extrapolation step sizes. Only a small transient error appears as $\Delta t$ increases. This is already an efficiency gain compared to the microscopic simulation: we obtain the same steady-state value while taking larger time steps. The right plot of Figure \ref{fig:triatommicro-macro} contains a surprising result. Micro-macro acceleration with extrapolating only the mean of $\xi_2$ completely eliminates the modeling error made by the effective dynamics of $\xi_2$, even though the mean of $\xi_2$ contains some fast dynamics. This is the case for multiple extrapolation step sizes. One reason why micro-macro acceleration is more accurate than the effective dynamics is that the microscopic dynamics is taken into account as well during the simulation stage. A complete explanation of the high level of accuracy for large $\Delta t$ is however not yet available.

We can conclude, numerically, that micro-macro acceleration yields accurate simulations results when extrapolating either one of the two reaction coordinates. Micro-macro acceleration attains the same level of accuracy or is more accurate than the effective dynamics. At the same time, we are able to bridge the gap in the time-scale separation by taking extrapolation time steps that are about ten times larger than the microscopic time step. The time-scale separation is approximately a factor of ten between the slow and ten fast dynamics, due to the particular choice of values for $\varepsilon$ and $k$. We thus gain in efficiency compared to both the effective dynamics and the microscopic time integrator.

\subsubsection{A hierarchy of moments of $x_c$ and $y_c$} \label{subsubsec:badmac}
\begin{figure}
\centering
\begin{subfigure}[b]{0.5\textwidth}
\centering
\includegraphics[width=0.99\textwidth]{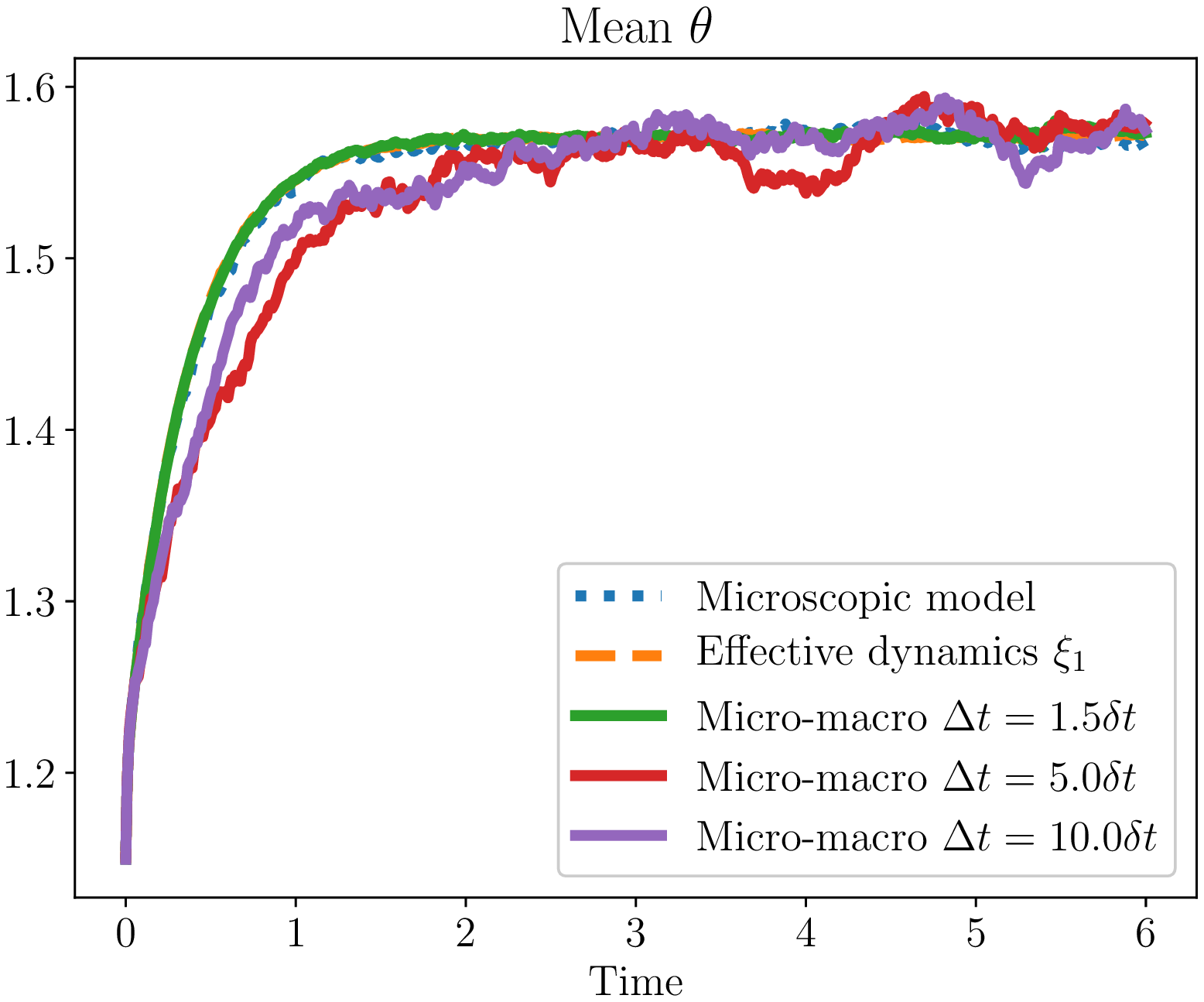}
\end{subfigure}%
\begin{subfigure}[b]{0.5\textwidth}
\centering
\includegraphics[width=0.99\textwidth]{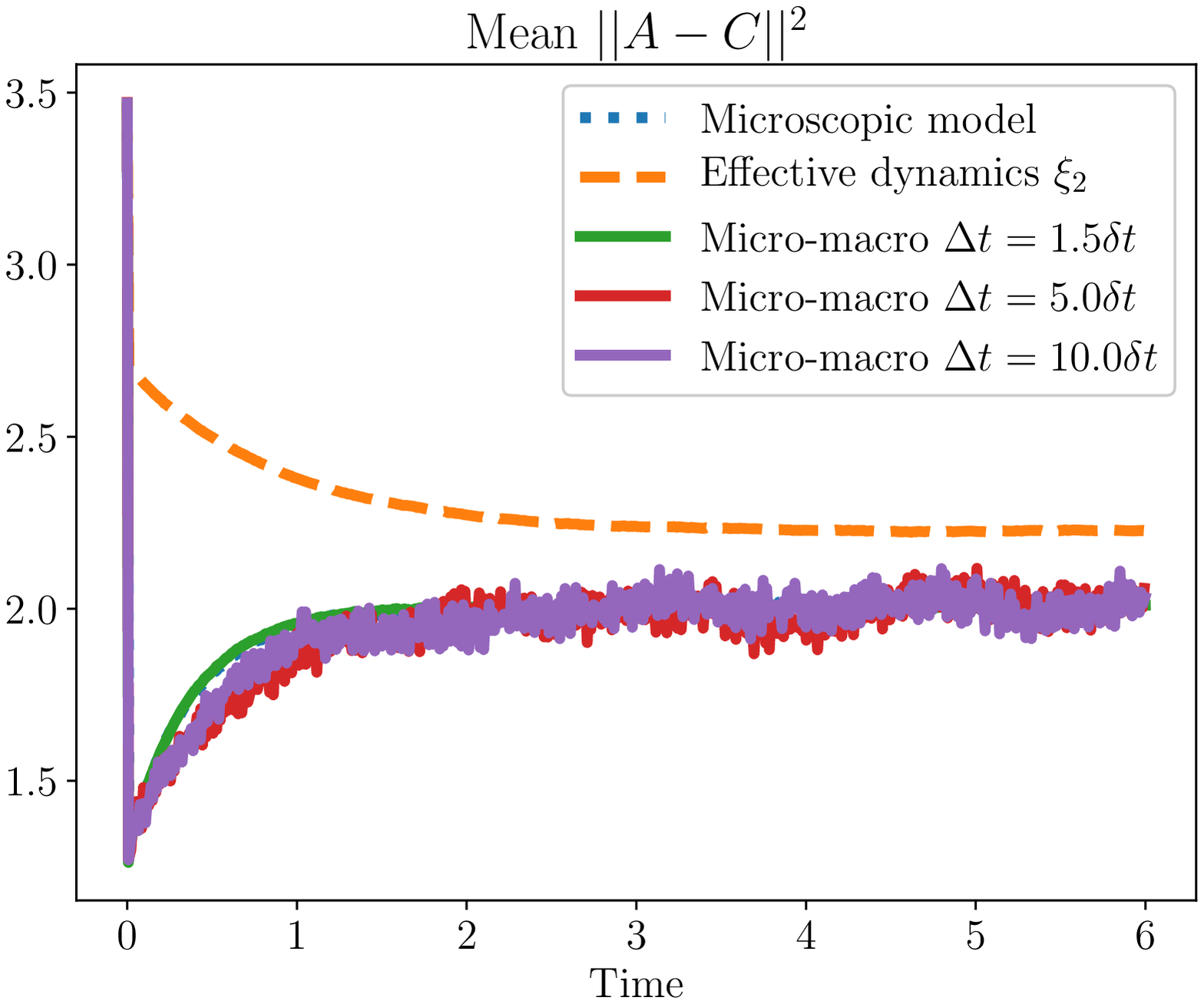}
\end{subfigure}
\caption{The evolution of the mean of $\theta$ (left) and the mean of $\norm{A-C}^2$ (right) computed with micro-macro acceleration and extrapolating the first two moments of $x_c$ and $y_c$ as macroscopic states. The noise amplitude on $\xi_1$ and $\xi_2$ increases when $\Delta t$ increases.}
\label{fig:triatom_moments}
\end{figure}

Reaction coordinates $\xi_1$ and $\xi_2$ are very closely connected (or identical to) the angle $\theta$, and micro-macro acceleration yields very accurate results when extrapolating the mean of either reaction coordinates. We now consider a set of pure moments of the three-atom molecule~\eqref{eq:triatommolecule} that also determines the angle $\theta$ uniquely and investigate the numerical results. Suppose we take the first and second moments of $x_c$ and $y_c$ as macroscopic state variables
\[
\mathbf{R}(\bm{X}) = \left(x_c, y_c, x_c^2, y_c^2\right).
\]
Figure \ref{fig:triatom_moments} depicts the numerical approximations by micro-macro acceleration with these moments as macroscopic state variables, with the same large extrapolation time steps and $N=5 \cdot 10^4$ particles.

The numerical results indicate that micro-macro acceleration follows the exact microscopic dynamics of $\theta$ and $\norm{A-C}^2$ well, for small extrapolation steps, i.e., $\Delta t = 1.5\delta t$. However, when the extrapolation step grows, micro-macro acceleration remains unbiased but the statistical noise increases. The noise on $\xi_2$ in the right panel of Figure \ref{fig:triatom_moments} is higher than the noise on $\xi_1$. This effect is mainly due to the fact that we are not extrapolating any states of $x_a$, such that $x_a$ can move freely. In the experiments in previous section, both $\xi_1$ and $\xi_2$ depend on $x_a$ and the noise is a lot lower. We leave the quantification of statistical error through micro-macro acceleration for future work.

\subsection{Conclusion on the number of macroscopic state variables} \label{subsec:conlcusomoments}
For both FENE-dumbbells and the three-atom molecule, we can conclude that choosing a set of macroscopic state variables that is in close connection to some quantity of interest (the angle $\theta$ for the three-atom molecule and the stress tensor $\tau$ for FENE-dumbbells) yields accurate simulation results. In case of FENE-dumbbells, the third hierarchy of macroscopic states approximates the stress very well, with fewer moments than the other hierarchies of macroscopic state variables. Also, in case of the three-atom molecule, micro-macro acceleration can improve on the effective dynamics for reaction coordinates that are related to the quantity of interest. 

As a second observation, we need to ensure to extrapolate enough macroscopic state variables related to the quantity of interest to keep the statistical error low enough. When extrapolating only moments of atom $C$ in the three-atom molecule, the noise on the angle $\theta$ is smaller than on the distance $\norm{A-C}^2$ since the variable $x_a$ can move around freely.

\section{Impact of the extrapolation step size on accuracy} \label{sec:extraoplation}
In this Section, we look in more detail into the effect of the extrapolation time step $\Delta t$ on accuracy. We assume that the macroscopic state variables are given, for example through an analysis similar to those in Section ~\ref{sec:states}. In the following examples, we are interested in the maximal extrapolation step $\Delta t$ we can employ, such that the micro-macro acceleration error remains within a given tolerance, specifically as a function of the time-scale separation present in the system. For systems in which the time-scale separation can be expressed in terms of a small-scale parameter $\varepsilon$, an approximate macroscopic model usually exists for the slow dynamics \cite{pavliotis2008multiscale}, that converges to the microscopic dynamics when $\varepsilon$ to 0. The goal of this Section is to determine the maximal $\Delta t$ such that micro-macro acceleration is more accurate than the approximate macroscopic model. In Section~\ref{subsec:bimodal}, we first investigate the efficiency gain on a slow-fast system with a double-well potential, and in Section~\ref{subsec:periodic}, we discuss a simple linear SDE with a periodic external force.

\subsection{A slow-fast bimodal system} \label{subsec:bimodal}
\begin{figure}
\centering
\includegraphics[width=0.4\textwidth]{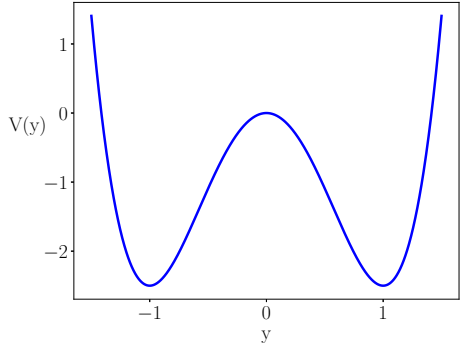}
\caption{The double well potential $V(y) = \frac{1}{4\varepsilon}y^4-\frac{1}{2\varepsilon} y^2$ for the slow-fast bimodal system \eqref{eq:bimodalsde} with $\varepsilon = 0.1$.}
\label{fig:doublewell}
\end{figure}

As the first example, we consider a slow-fast system, where the fast component $Y(t)$ is governed by a double-well potential $V(y) = \frac{1}{4 \varepsilon}y^4 - \frac{1}{2 \varepsilon}y^2$
\begin{equation} \label{eq:bimodalsde}
\begin{aligned}
dX &= -(2X+Y)dt + 0.1dW_x \\
dY &= \frac{1}{\varepsilon}(Y-Y^3)dt + \frac{1}{\sqrt{\varepsilon}} dW_y = -\frac{1}{\varepsilon}\nabla V(Y)dt + \frac{1}{\sqrt{\varepsilon}} dW_y.
\end{aligned}
\end{equation}
The potential energy function $V(y)$ is depicted in Figure \ref{fig:doublewell}.

\begin{figure}
\centering
\begin{subfigure}[b]{0.5\textwidth}
\centering
\includegraphics[width=0.8\textwidth]{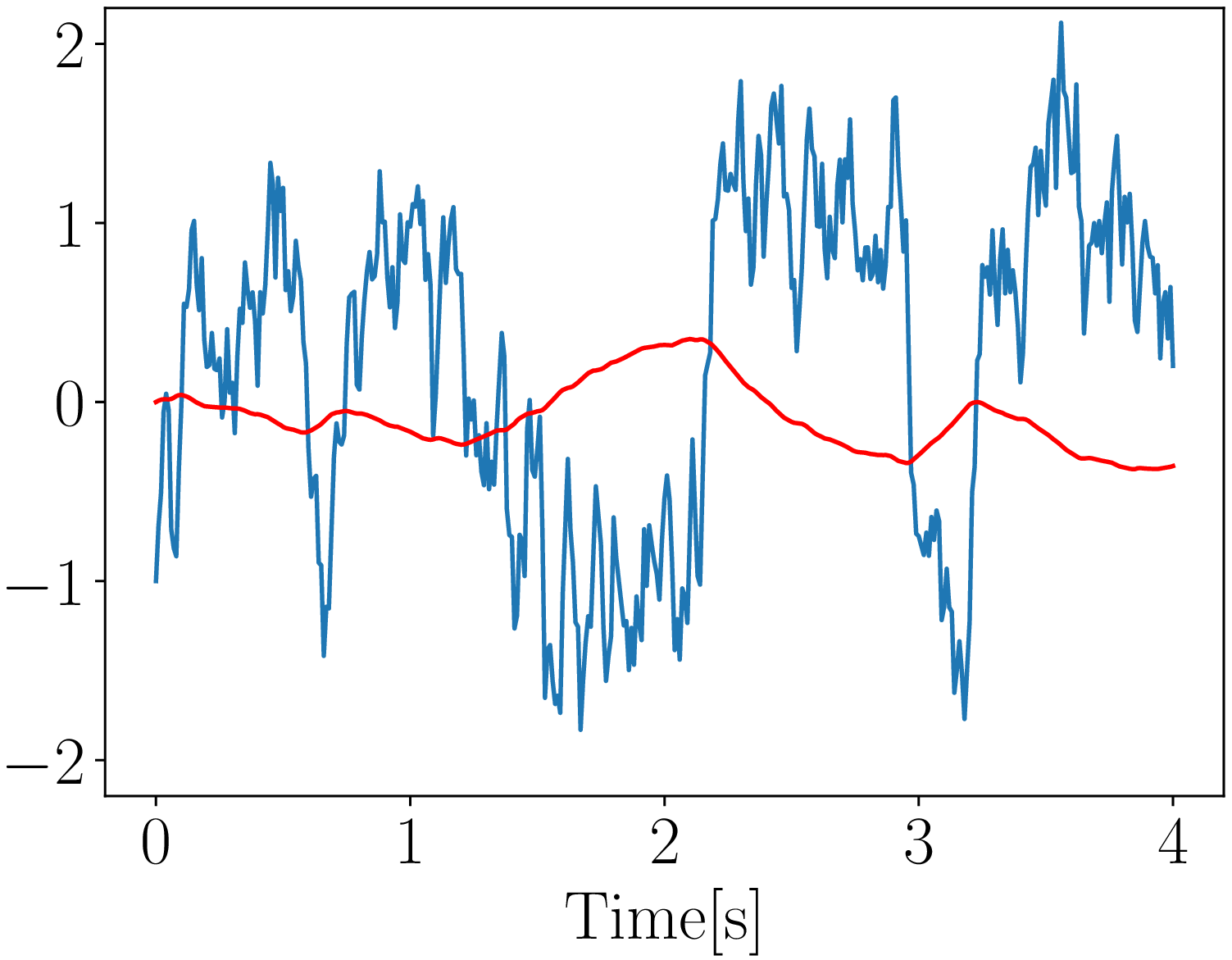}
\end{subfigure}%
\begin{subfigure}[b]{0.5\textwidth}
\centering
\includegraphics[width=0.83\textwidth]{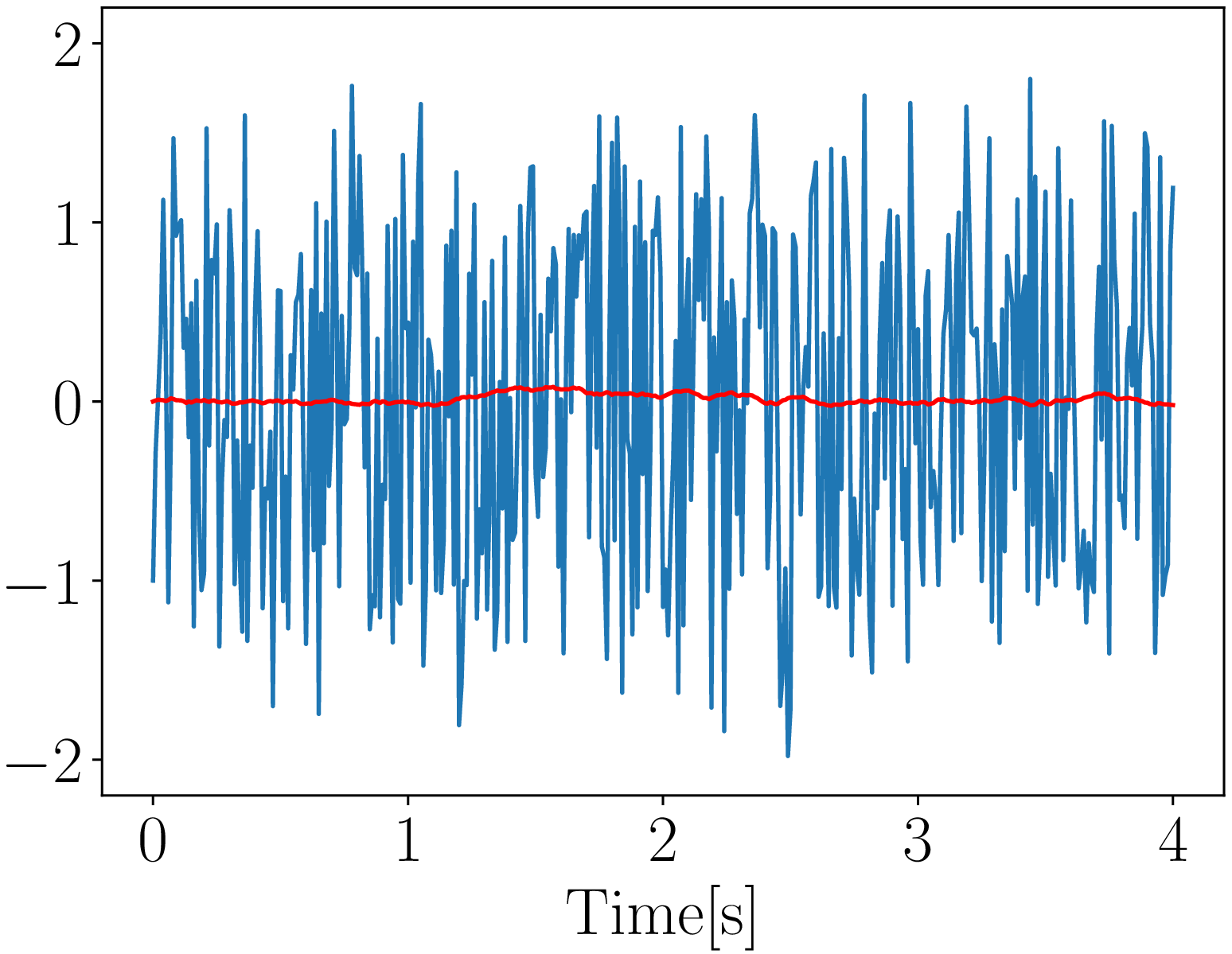}
\end{subfigure}
\caption{The motion of one particle through the bimodal system \eqref{eq:bimodalsde}. For large $\varepsilon$ (left), the fast component (blue) influences greatly the slow component (red) due to longer residence in the wells. When $\varepsilon$ is small (right), the slow component remains almost unaffected by the fast switching of the fast mode.}
\label{fig:particlemovement}
\end{figure}

With this double-well potential, the fast component will reside in one well for some time and then switch at a random instant in time to the other well due to the Brownian motion term in the evolution equation for $Y(t)$. The mean switching time between the two wells depends on $\varepsilon$ and is derived in detail in \cite{bruna2014model}. Figure \ref{fig:particlemovement} depicts the switching process of one particle for a small $\varepsilon = 0.001$ (left) and a large $\varepsilon = 0.1$ (right). For large time-scale separations (small $\varepsilon$), the fast component switches frequently between the wells, so that the slow component (in red) is not heavily influenced by the switching. The slow component only feels the invariant distribution of the position of the fast component. For low time-scale separations, the fast component resides longer in one well before switching to the other. In this case, the slow component feels in which well the fast component resides and evolves accordingly.

The slow-fast bimodal system in this section differs from the bimodal behavior of the three-atom molecule. In the current example, the bimodal motion is fast while in the previous example it is slow. Also, the effective dynamics from the three-atom molecule \eqref{eq:effecitvedynamics} is not applicable here since there is an $\varepsilon$ present in the diffusion term of \eqref{eq:bimodalsde}. 

In Section~\ref{subsec:pavliotis}, we introduce an approximate macroscopic model for the slow component of general slow-fast systems, based on averaging out the fast component. We show the accuracy of the approximate macroscopic model for small and larger values of $\varepsilon$. Finally, in Section~\ref{subsec:bimodalmM}, we discuss the accuracy and efficiency gain of micro-macro acceleration compared to the approximate macroscopic model, applied to the bimodal system~\eqref{eq:bimodalsde}.

\subsubsection{An approximate macroscopic model for general slow-fast systems} \label{subsec:pavliotis}
For general slow-fast systems, one can derive an approximate macroscopic model for the slow component by averaging out the mean of the fast dynamics $Y$, given value of the slow component $X$. We then substitute the averaged value of $Y$ into the evolution equation for the slow component. In case of the bimodal system, the fast dynamics is autonomous and the invariant measure is \cite{pavliotis2008multiscale}
\[
\mu_{\infty}(y | x) = \frac{1}{Z} \exp\left( - 2 \varepsilon V(y)\right).
\]
The mean of the $y$, given a value of $x$ equals zero
\[
\int_{\mathbb{R}} y \ \mu_{\infty}(y | x) dy = 0,
\]
so that the approximate macroscopic model for the bimodal system \eqref{eq:bimodalsde} reads
\begin{equation} \label{eq:approximatebimodal}
d\overset{\_}{X} = -2\overset{\_}{X} dt + 0.1 dW_x.
\end{equation}
Equation \eqref{eq:approximatebimodal} converges weakly to \eqref{eq:bimodalsde} when $\varepsilon$ decreases to 0 \cite{pavliotis2008multiscale}, implying that the approximate macroscopic model makes a modeling error whenever $\varepsilon \neq 0$. The reason for a modeling error can easily be identified from the left plot of Figure \ref{fig:particlemovement}. Indeed, as we mentioned, for larger $\varepsilon$ the fast mode heavily influences the slow mode, so that we cannot just replace the fast behavior by its mean. 

\begin{figure}
\begin{subfigure}[b]{0.5\textwidth}
\centering
\includegraphics[width=0.9\textwidth]{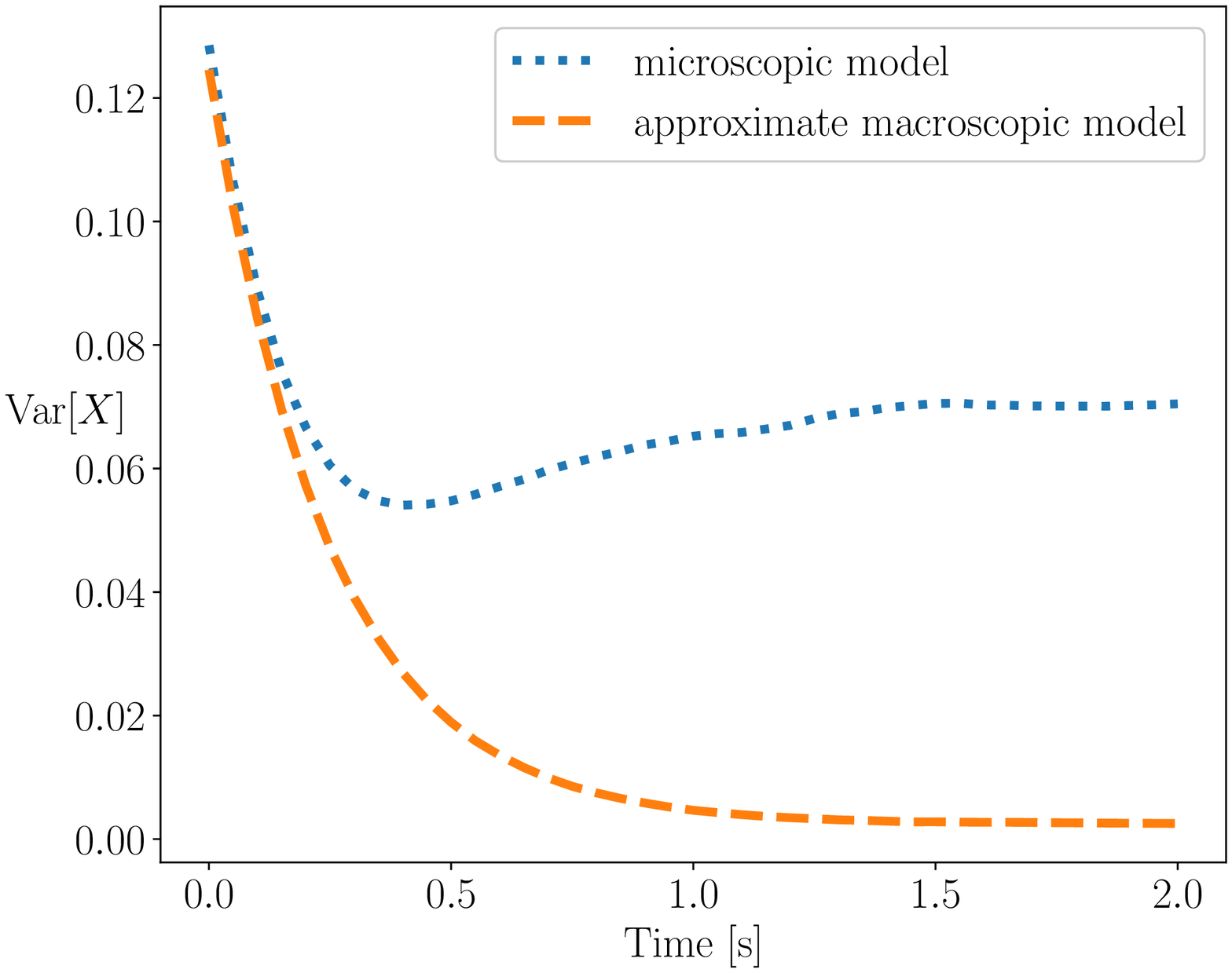}
\end{subfigure}%
\begin{subfigure}[b]{0.5\textwidth}
\centering
\includegraphics[width=0.9\textwidth]{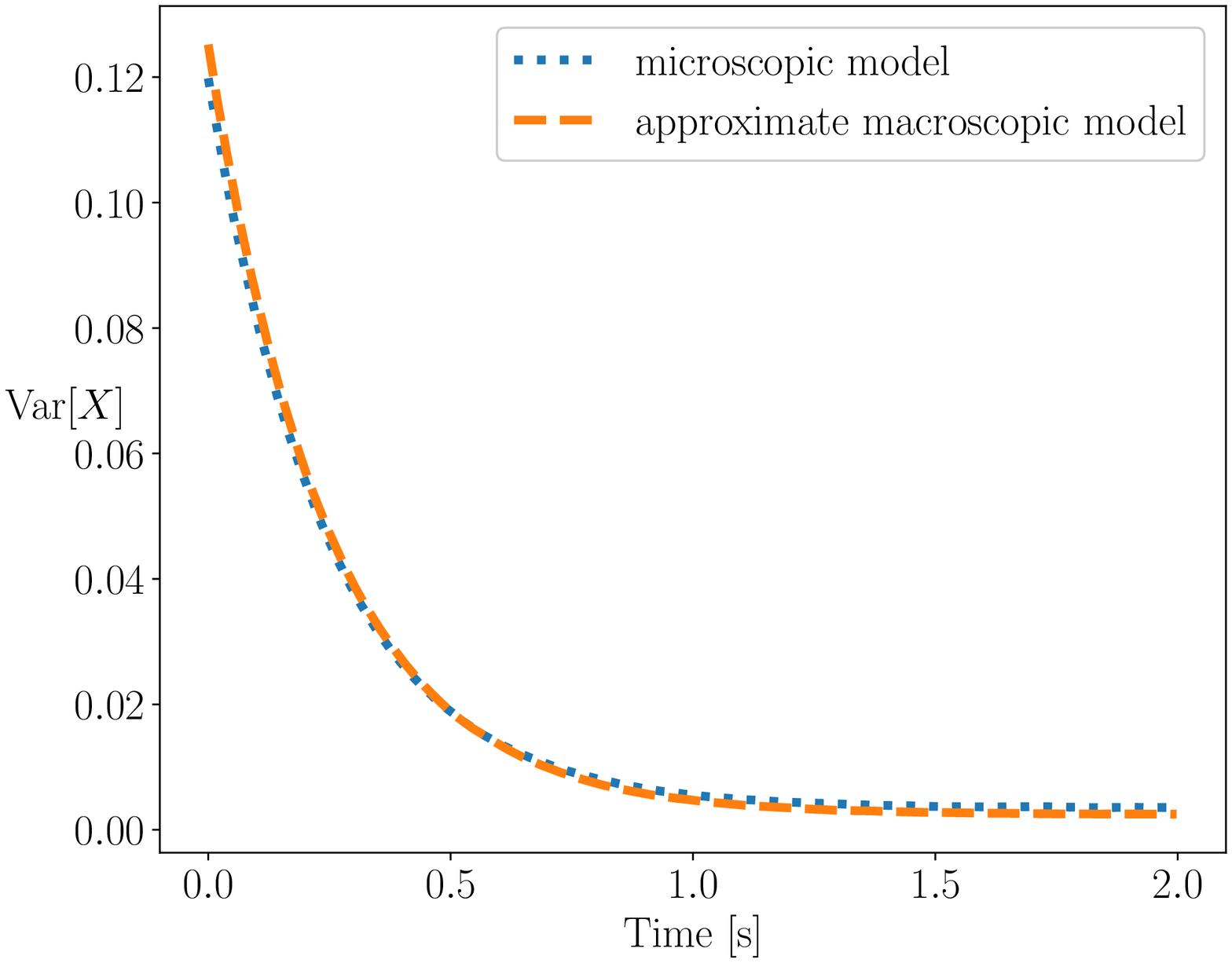}
\end{subfigure}
\caption{Variance of the slow component $X$ in \eqref{eq:bimodalsde}, computed with the full microscopic model(blue) and the approximate macroscopic model (orange), specifically for $\varepsilon=0.1$ (left) and $\varepsilon=0.001$ (right). For large $\varepsilon$, the approximate macroscopic model makes a large steady-state error compared to the microscopic dynamics. For small $\varepsilon$, the approximate macroscopic model is more accurate.}
\label{fig:bimodalmodelingerror}
\end{figure}

There is a large variance in $X$ induced by the switching of $Y$ that we ignore by substituting the mean of $Y$ for $Y$ itself. When $\varepsilon$ is lower, the error also lowers as we identified in the right plot of Figure \ref{fig:particlemovement}. This induced variance is illustrated in Figure \ref{fig:bimodalmodelingerror}, which shows the variance of the slow component of \eqref{eq:bimodalsde} as a function of time, for $\varepsilon = 0.1$(left) and $\varepsilon = 0.001$(right). The step size is for both models is $\delta t = \varepsilon/10$ and the initial condition is far from the equilibrium distribution of the system.

The left plot of Figure \ref{fig:bimodalmodelingerror} indicates a large steady-state error on the slow variance for moderate $\varepsilon$; this error is also present for $\varepsilon=0.001$ but a lot smaller. With the micro-macro acceleration method, we aim to reduce this steady-state error, especially for moderate time-scale separations.

\subsubsection{Efficiency gain by micro-macro acceleration} \label{subsec:bimodalmM}
\begin{figure}
\centering
\begin{subfigure}[b]{0.5\textwidth}
	\centering
\includegraphics[width=0.95\textwidth]{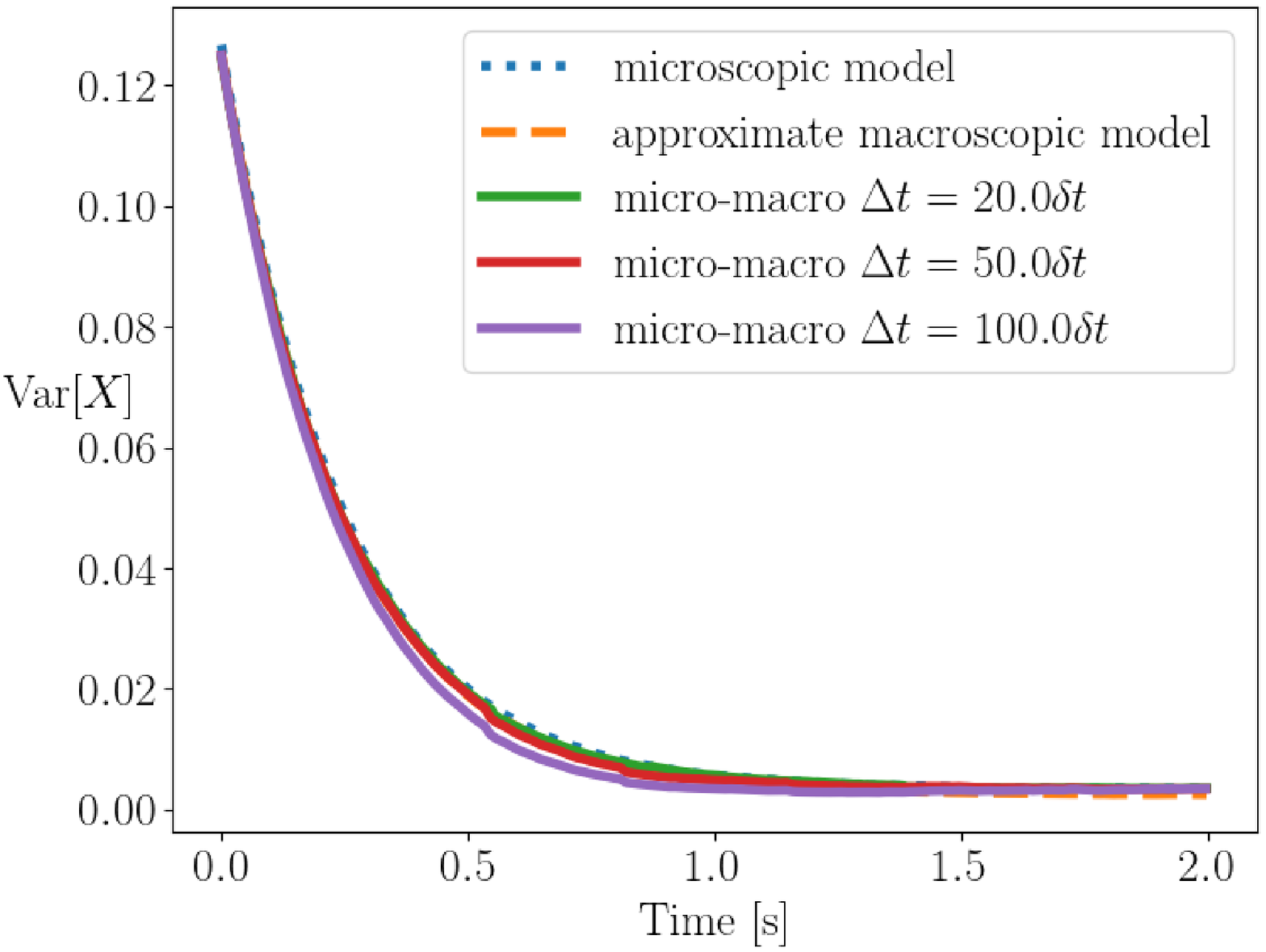}
\end{subfigure}%
\begin{subfigure}[b]{0.5\textwidth}
	\centering
	\includegraphics[width=0.95\textwidth]{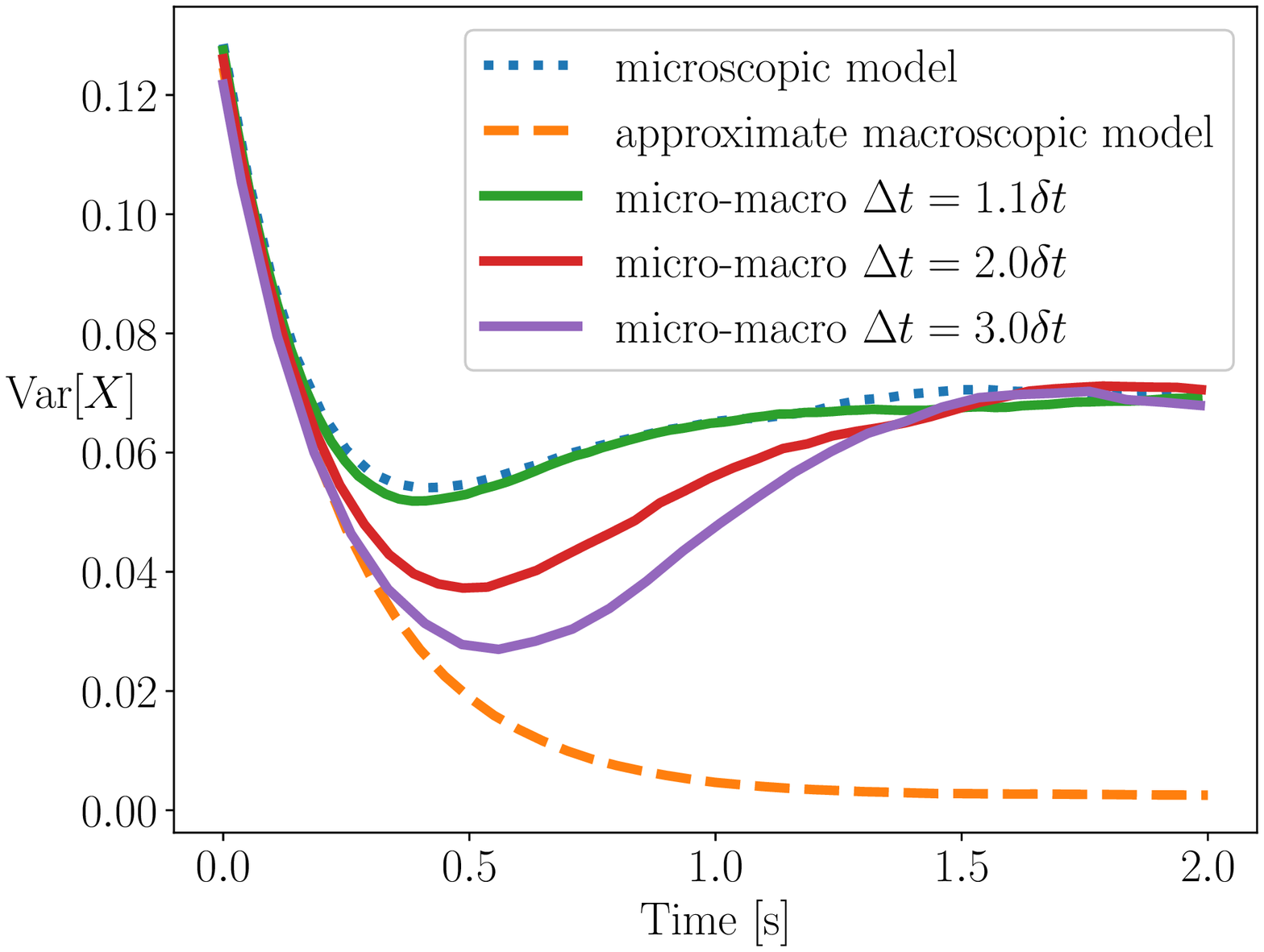}
\end{subfigure}%
\caption{Variance of the slow component of \eqref{eq:bimodalsde}, computed with the microscopic model (blue), approximate macroscopic model (orange) and micro-macro acceleration with multiple extrapolation steps (full lines) for $\varepsilon=0.001$ (left) and $\varepsilon=0.1$ (right). For a large time-scale separation, micro-macro acceleration follows the exact dynamics well, even for large extrapolation step sizes. When the time-scale separation is low, the possible gain is low and micro-macro acceleration reaches the correct steady-state value but the transient error increases as $\Delta t$ increases.}
\label{fig:variance_mM}
\end{figure}

In the following experiments, we run the micro-macro acceleration method on the bimodal system \eqref{eq:bimodalsde}, and extrapolate the first and second moment of $X(t)$, with $N=10^5$ particles. Both moments are required, since we need to approximate the variance of the slow component as closely as possible. We consider values for the time-scale parameter: $\varepsilon=0.001$ and $\varepsilon = 0.1$. In the former case, we expect a large efficiency gain due to a high time-scale separation. In the latter case, the time-scale separation is only moderate and the expected efficiency gain is rather small. The numerical results are summarized in Figure \ref{fig:variance_mM}. On the one hand, for small $\varepsilon$, micro-macro acceleration can take extrapolation steps that are a hundred times larger than the microscopic step size, while only creating a small transient error. On the other hand, when $\varepsilon$ is large, micro-macro acceleration is able to eliminate the error of the approximate macroscopic model and reach the same steady state value for the variance as the microscopic integrator. At the same time, micro-macro acceleration also gains in efficiency compared to the microscopic integrator by taking larger time steps. There, however, remains a transient error on the variance of the slow mode, that enlarges as $\Delta t$ increases. There is currently no error control available for micro-macro acceleration and we leave this topic for further research.

To conclude, with micro-macro acceleration, we can almost bridge the time-scale separation gap by taking extrapolation steps on the order of the small dynamics then $\varepsilon$ is small, and when $\varepsilon$ is large, we can eliminate the steady-state error of the approximate macroscopic model with time steps that are slightly larger than those of the microscopic time integrator. Thus, on this double-well system, we gain in efficiency for all values of $\varepsilon$.

\subsection{Accuracy analysis of extrapolation on a periodically driven linear system}  \label{subsec:periodic}
For the final numerical example of this manuscript, we consider a slow-fast linear system where we add a periodic external force to evolution equation of the the slow component \cite{li2008effectiveness},
\begin{equation} \label{eq:lineardrivensde}
\begin{cases}
dX = -2(X+Y)dt + 10\sin(2\pi t)dt + dW_x \\
dY = \frac{1}{\varepsilon}(X-Y)dt + \frac{1}{\sqrt{\varepsilon}}dW_y.
\end{cases}
\end{equation}
Hence, the mean of $X$ has a periodic invariant solution. The periodic external force makes it easier to compare the numerical solution of different methods to the exact solution. The errors between different models are present in every period of the invariant solution and do not damp out to zero as time increases. We can thus compute errors of different methods by simply computing $L_2$-norm of their difference with the exact solution, over one period of the solution.

We can also apply the averaging strategy from the previous example \eqref{eq:bimodalsde} to the periodic slow-fast model \eqref{eq:lineardrivensde}. The invariant distribution of $Y$, given a value for $X$ is Gaussian and reads
\[
\mu_{\infty}(y | x) = \frac{1}{Z} \exp\left(2(y-x)^2\right),
\]
with $Z$ the normalization constant. The mean of $Y$, computed against its invariant distribution, for a given value of $X$, thus reads
\[
\int_{\mathbb{R}} y \ \mu_{\infty}(y | x) dy = x.
\]
Therefore, the approximate macroscopic model for the slow component $X(t)$ is
\begin{equation} \label{eq:linearapproximate}
d\overset{\_}{X} = -4\overset{\_}{X}dt +10\sin(2\pi t)dt+ dW_x.
\end{equation}

\begin{figure}
	\begin{subfigure}[b]{0.5\textwidth}
		\centering
		\includegraphics[width=0.95\textwidth]{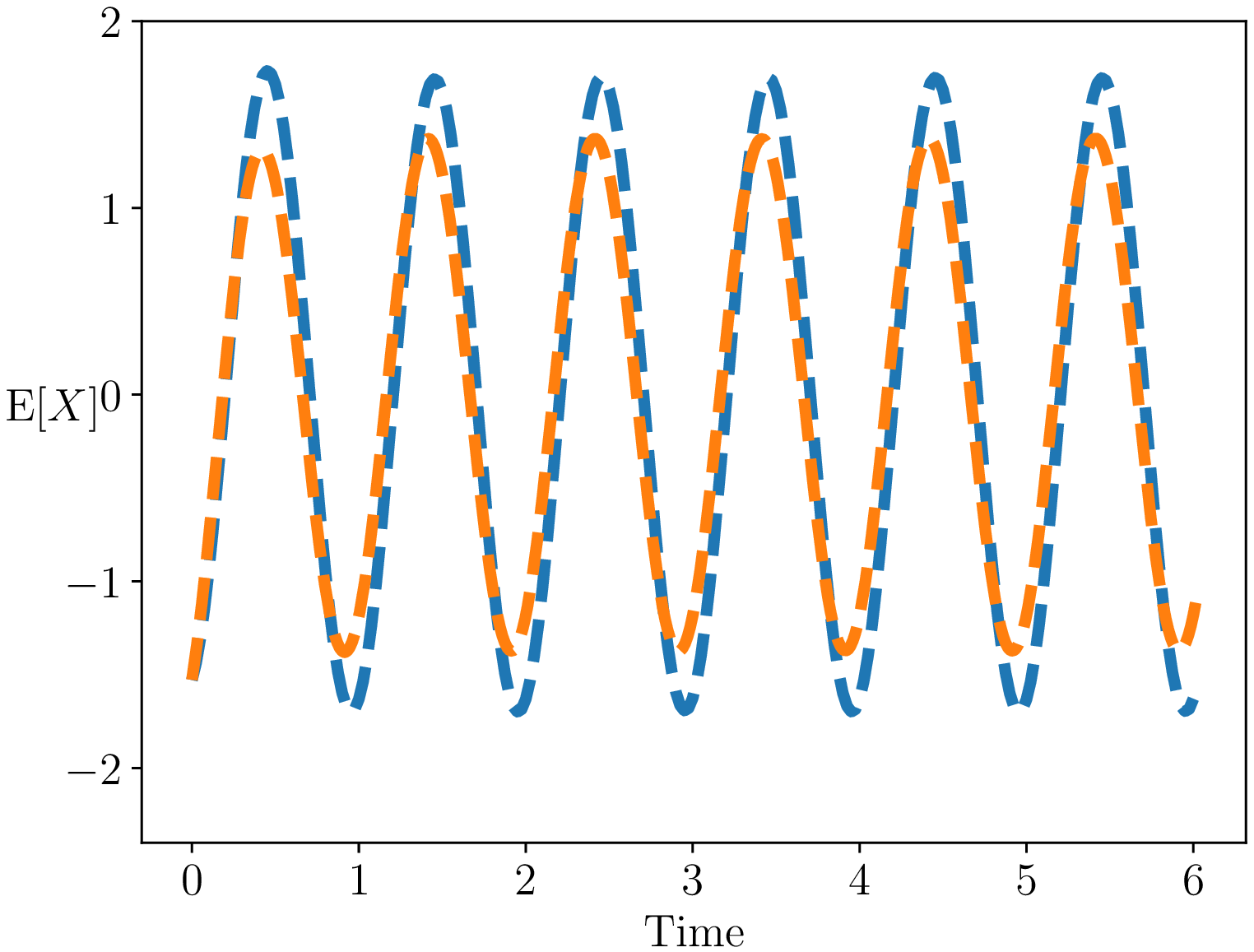}
	\end{subfigure}%
	\begin{subfigure}[b]{0.5\textwidth}
		\centering
		\includegraphics[width=0.95\textwidth]{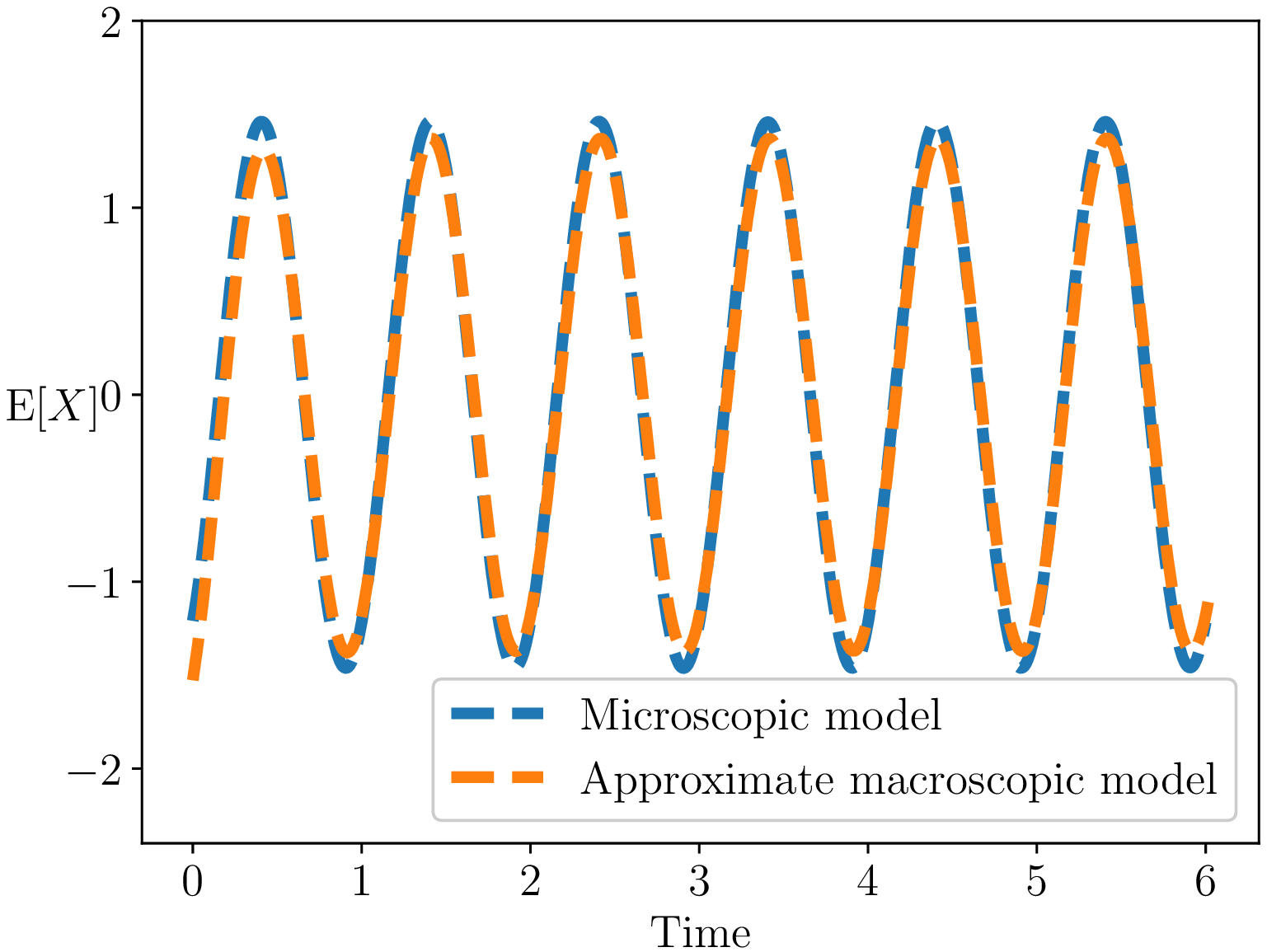}
	\end{subfigure}
	\caption{The evolution of the mean of $X(t)$ in \eqref{eq:lineardrivensde} (blue) and of the approximate macroscopic model in \eqref{eq:linearapproximate} (orange). The values for the time-scale separation are $\varepsilon=0.5$ (left) and $\varepsilon=0.05$ (right).}
	\label{fig:periodicmodelingerror}
\end{figure}

\begin{figure}
    \centering
    \includegraphics[width=0.45\textwidth]{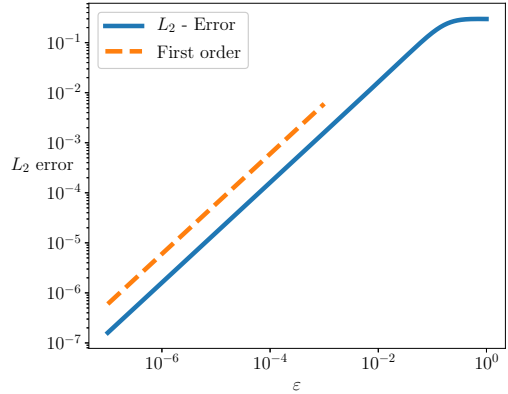}
    \caption{The $L_2$ error made by the approximate macroscopic model~\eqref{eq:linearapproximate} (in blue) decreases linearly with $\varepsilon$.}
    \label{fig:linear_error}
\end{figure}

Let us first compare the numerical solution of the approximate macroscopic model~\eqref{eq:linearapproximate} to the exact microscopic dynamics~\eqref{eq:lineardrivensde}, for different values of $\varepsilon$. For all the numerics in this Section, the initial condition is taken on the invariant curve of the exact continuous solution, the microscopic time step is $\delta t = \varepsilon / 10$ and we use $N=10^5$ Monte Carlo particles. Figure \ref{fig:periodicmodelingerror} depicts the evolution of the mean of $X$ for $\varepsilon=0.5$ on the left, and $\varepsilon=0.05$ on the right. The numerical results show that the approximate macroscopic model (in blue) makes an error, compared to the exact microscopic model (in orange). The error is mostly visible in the amplitude of the periodic solution, and the discrepancy decreases when $\varepsilon$ decreases to zero. We show in Appendix~\ref{app:B} that the error between the microscopic and approximate macroscopic model decreases linearly with $\varepsilon$, as is also visible on Figure~\ref{fig:linear_error}.

\subsubsection{Tuning the error of micro-macro acceleration}
Let us now investigate the accuracy of micro-macro acceleration for the periodically driven linear system~\eqref{eq:lineardrivensde}. For every value of $\varepsilon$, we want to examine the maximal extrapolation step micro-macro acceleration can take, while keeping the error smaller than the corresponding error of the approximate macroscopic model. 

\begin{figure}
\centering
\begin{subfigure}[b]{0.43\textwidth}
\includegraphics[width=\textwidth]{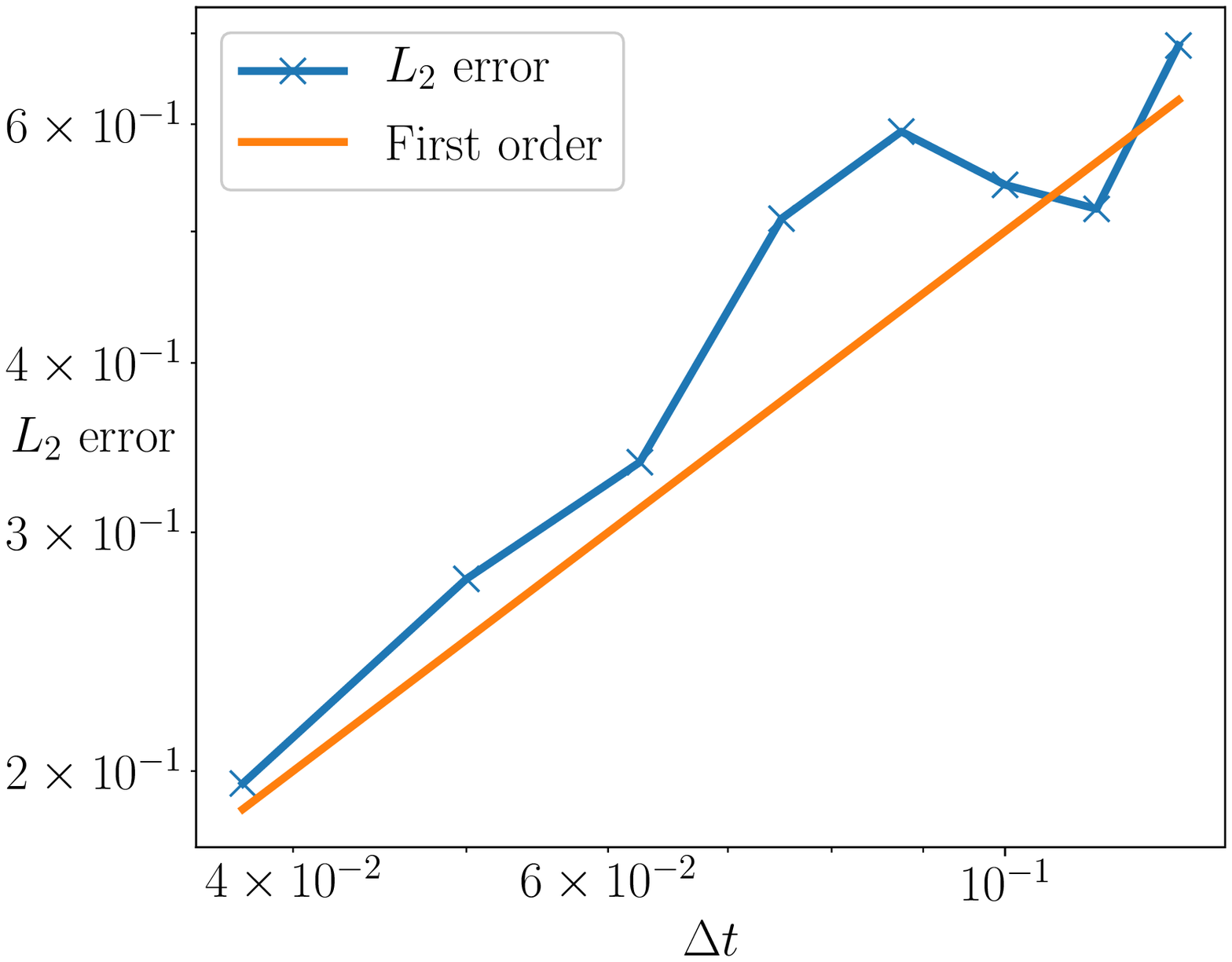}
\end{subfigure}%
\begin{subfigure}[b]{0.42\textwidth}
\includegraphics[width=\textwidth]{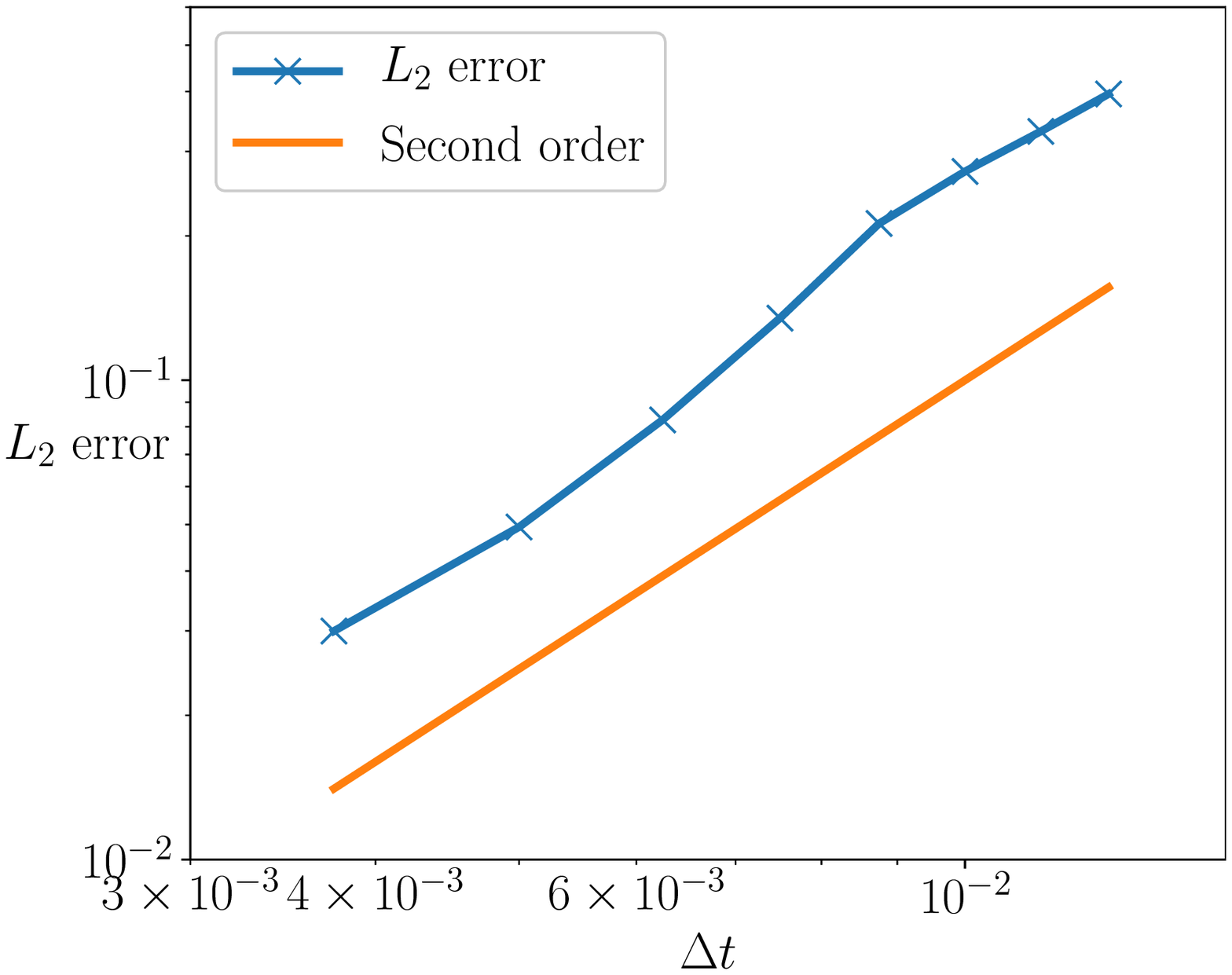}
\end{subfigure}
\caption{Convergence of micro-macro acceleration as $\Delta t$ decreases for $\varepsilon=0.5$ (left) and $\varepsilon=0.05$ (right). For large $\varepsilon$ the convergence is first order, while it is second for small $\varepsilon$.}
\label{fig:convergencemM}
\end{figure}

To determine the maximal extrapolation step, we first need to know how the micro-macro acceleration error decreases as $\Delta t$ decreases to $\delta t$. As of now, there exist no theoretical result on the convergence rate of micro-macro acceleration with relative entropy matching \cite{lelievre2018analysis}, so we can only investigate the rate of convergence numerically. In Figure \ref{fig:convergencemM}, we show the error of micro-macro acceleration as a function of the extrapolation time step $\Delta t$, for the same values of $\varepsilon$ as in Figure \ref{fig:periodicmodelingerror}. For large values of $\varepsilon$, the micro-macro acceleration error decreases linearly with $\Delta t$ as is shown on the left plot in Figure \ref{fig:convergencemM}. We expect first-order convergence since the linear extrapolation of macroscopic state variables mimics the forward Euler method. However, for small $\varepsilon$, we observe that the error decreases quadratically with $\Delta t$.

Having the convergence order of both the approximate macroscopic model and micro-macro acceleration for small $\varepsilon$ in place, we will derive an expression for the maximal extrapolation step we can take,  before the approximate macroscopic model becomes more accurate than micro-macro acceleration. Using the linear convergence of the approximate macroscopic model, and the quadratic convergence of micro-macro acceleration, both methods attain the same accuracy when
\begin{equation} \label{eq:crossovereq}
\alpha(\varepsilon)(\Delta t_{\text{max}})^2 = \beta \varepsilon,
\end{equation}
with $\Delta t_{\text{max}}$ the maximal extrapolation step possible before the approximate macroscopic model is more accurate than micro-macro acceleration. We define $\beta$ as the constant of proportionality for the approximate macroscopic model and the function $\alpha(\varepsilon)$ as the proportionality constant for micro-macro acceleration, which can in principle depend on $\varepsilon$. For a more natural comparison of the maximal extrapolation step for different values of $\varepsilon$, we define the dimensionless `extrapolation factor' $M$ as $\Delta t =M\delta t$. Furthermore, due to the stability requirement for the inner integrator, we choose the inner microscopic step $\delta t$ proportional to $\varepsilon$,  $\delta t = C \varepsilon$. Putting all the parameters together, equation \eqref{eq:crossovereq} yields
\begin{equation} \label{eq:reducedcrossovereq}
\alpha(\varepsilon)C^2(M_{\text{max}} \varepsilon)^2 = \beta \varepsilon,
\end{equation}
such that a direct expression for $\Delta t_{\text{max}} = M_{\text{max}} \delta t$ follows. Unfortunately, we lack a direct expression for $\alpha(\varepsilon)$, so we cannot use \eqref{eq:reducedcrossovereq} directly. We must resort to numerical experiments to deduce the maximal extrapolation step $\Delta t_{\text{max}}$ and the maximal extrapolation factor $M_{\text{max}}$. 

\begin{figure}
	\centering
	\includegraphics[width=0.65\textwidth]{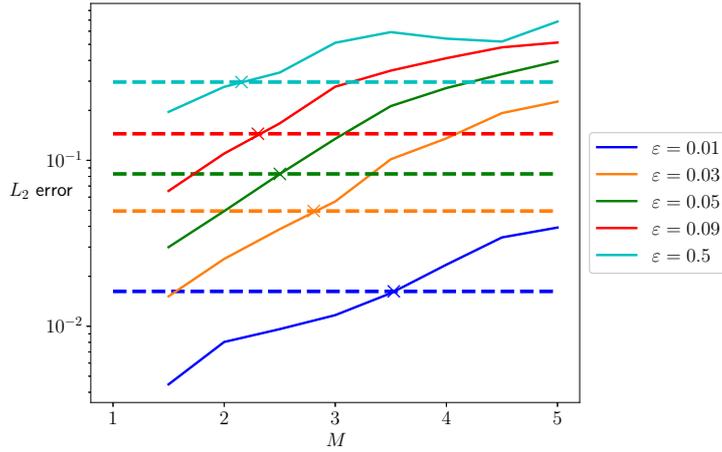}
	\caption{The dashed lines indicate the error of the approximate macroscopic model for the given value of $\varepsilon$. The solid line of the same color denotes the error of micro-macro acceleration as a function of the extrapolation factor $M$. The crosses indicate where both methods obtain the same accuracy for that value of $\varepsilon$.}
	\label{fig:crossover}
\end{figure}

\begin{figure}
	\centering
	\includegraphics[width=0.5\textwidth]{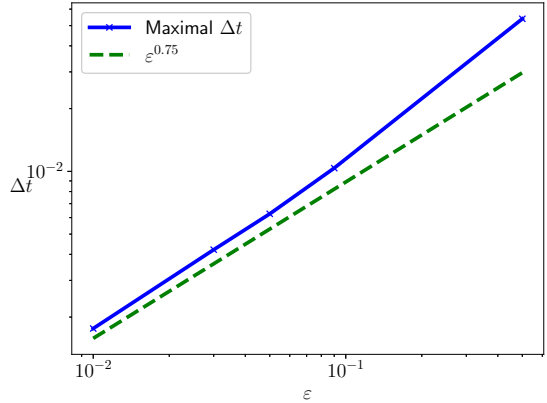}
	\caption{The maximal extrapolation step $\Delta t_{\text{max}}$, before the approximate macroscopic model becomes more accurate than micro-macro acceleration, decreases slower than linearly with $\varepsilon$.}
	\label{fig:deltatdecrease}
\end{figure}

In Figure \ref{fig:crossover}, we plot in solid lines the micro-macro acceleration error for multiple values of $\varepsilon$ as a function of the extrapolation factor $M$. The error of the approximate macroscopic model is also given by the dashed lines, and the crosses indicate where both errors are equal. Denote the associated maximal extrapolation factor by $M_{\text{max}}(\varepsilon)$. Whenever $M < M_{\text{max}}(\varepsilon)$, the micro-macro acceleration method is more accurate than the approximate macroscopic model for that value of $\varepsilon$, and vice versa when $M > M_{\text{max}}(\varepsilon)$. The numerical results indicate that the maximal extrapolation factor $M_{\text{max}}$ increases as $\varepsilon$ decreases.

An important consequence of the the fact that $M_{\text{max}}$ increases, is that the maximal extrapolation step $\Delta t_{\text{max}}$ decreases slower than linearly with $\varepsilon$, since the microscopic inner step $\delta t$ decreases linearly with $\varepsilon$. Indeed, Figure \ref{fig:deltatdecrease} depicts that $\Delta t_{\text{max}}$ decreases approximately as $\varepsilon^{0.75}$. As a result, micro-macro acceleration can be more accurate than the approximate macroscopic model if we choose $\Delta t = \mathcal{O} (\varepsilon)$. We can also tune the micro-macro acceleration error. Given a tolerance, or maximal value for the error, we can choose the extrapolation step $\Delta t$ such the error is beneath the tolerance, while still taking larger time steps than the microscopic time integrator. We conclude that micro-macro acceleration gains in efficiency compared to the approximate macroscopic model by lowering the error, but also compared to the inner microscopic time integrator by taking larger time steps. The gain compared to the microscopic time integrator even becomes infinite when $\varepsilon$ decreases to zero.

\section{Conclusion} \label{sec:conclusion}
We studied the efficiency of a micro-macro acceleration procedure. We introduced the micro-macro acceleration algorithm in mathematical detail in Section~\ref{sec:mM} and discussed some implementation issues in Section~\ref{sec:impl_det}.

\vspace{2mm}
We investigated the effect of the number and nature of the macroscopic state variables in detail in Section~\ref{sec:states}. The state variables have a big impact on the efficiency of micro-macro acceleration. Increasing the number of state variables yields a higher accuracy of the matched distributions but increases the cost of matching.

The numerical results of Section~\ref{sec:states} show that a hierarchy of macroscopic state variables that is closely connected to some quantity of interest, gives accurate simulation results for that quantity. In case of FENE-dumbbells, the hierarchy of state variables that includes terms in the It{\^o} expansion of the stress tensor, gives better approximations to the stress tensor than other hierarchies. In fact, this hierarchy of states requires only four macroscopic state variables, while the others do not approximate the stress tensor, even with five macroscopic state variables. 

A similar conclusion holds for the three-atom molecule, where we worked with two reaction coordinates: the angle and the distance between the two outer atoms. Micro-macro acceleration with extrapolating the mean of the two reaction coordinates, follows the exact microscopic dynamics very well, even for quite large extrapolation step sizes. Moreover, the second reaction coordinates eliminates the steady-state error made by the approximate model based on the effective dynamics for that reaction coordinate. The reaction coordinate does not necessarily need to be slow, as we have seen with the second reaction coordinate, it just needs to be linked to the quantity of interest, the angle in this example. A hierarchy based only on moments also produced an unbiased results but the statistical error is larger. Thus, micro-macro acceleration achieves the same or a better accuracy than the effective dynamics, while allowing for larger time steps than the stiff microscopic solver, gaining efficiency compared to both models.

\vspace{2mm}
Besides the choice of macroscopic state variables, the extrapolation step size greatly influences the accuracy of micro-macro acceleration, as we studied in Section~\ref{sec:extraoplation}. For the slow-fast bimodal model, the approximate macroscopic model makes a modeling error on the variance of the slow component, especially visible when there is almost no time-scale separation. Micro-macro acceleration attains the correct steady-state value for different extrapolation step sizes, at the expense of a controllable transient error. We thus gain in accuracy compared to the approximate macroscopic model while gaining in efficiency against the microscopic integrator by taking larger time steps. For the slow-fast periodic model, we proved numerically that the efficiency gain by micro-macro acceleration becomes infinite as the time-scale separation increases to infinity. The maximal extrapolation step micro-macro acceleration can take before the approximate macroscopic model becomes more accurate decreases approximately as $\varepsilon^{0.75}$, while the microscopic time step must decrease linearly due to stability. Furthermore, we can also tune the extrapolation step such that the error stays below a given tolerance, while still taking larger time steps than the microscopic time integrator. The efficiency gain is thus also present in the periodic linear system.

\bibliographystyle{plain}
\bibliography{sample}

\appendix
\section{Derivation of the macroscopic state variables for FENE-dumbbells} \label{app:fene}
In this Section, we give a short derivation of the macroscopic state variables of the third hierarchy of the FENE-dumbbells process~\eqref{eq:strategy3}. We start by writing out the evolution equation of the first even moment $\mathbb{E}[X_t^2]$ and consequently add all terms that appear in the evolution equation. For all these new terms we also write the evolution equations and select extra terms and so forth.

First note that for any twice-differentiable function $f(x)$, the evolution of the expectation $\mathbb{E}[f(X_t)]$ is given by Itô's law
\begin{equation} \label{eq:evolutionfene}
\frac{d}{dt}\mathbb{E}[f(X_t)] = \kappa(t)\mathbb{E}[X_t f'(X_t)] - \frac{1}{2W_{\text{e}}} \mathbb{E}[F(X_t) f'(X_t)] + \frac{1}{2W_{\text{e}}} \mathbb{E}[f''(X_t)].
\end{equation}
Starting from the first even moment of the diffusion process $f_1(x) = x^2$, we find that
\[
\frac{d}{dt}\mathbb{E}[X_t^2] = 2\kappa(t)\mathbb{E}[X_t^2] - \frac{1}{W_{\text{e}}} \mathbb{E}\left[\frac{b^2 X_t^2}{b^2-X_t^2}\right] + \frac{1}{W_{\text{e}}}.
\]
Since the state variable $\mathbb{E}[X_t^2]$ appears in its evolution equation, we do not need to add it to the list of macroscopic state variables. We also do not need to add the constant term, so we add $\mathbb{E}\left[\frac{b^2 X_t^2}{b^2-X_t^2}\right]$ as the second macroscopic state variable for the FENE-process.

Let us now consider the second function $f_2(x) = \frac{x^2}{b^2-x^2}$. The first and second derivative are
\begin{align*}
\begin{split}
    f_2'(x) &= \frac{2xb^2}{(b^2-x^2)^2} \\
    f_2''(x) &= \frac{2b^2(b^2+3x^2)}{(b^2-x^2)^3}.
\end{split}
\end{align*}

Using equation \eqref{eq:evolutionfene}, the evolution equation of $\mathbb{E}[f_2(X_t)]$ hence becomes
\[
\frac{d}{dt}\mathbb{E}[f_2(X_t)] = \kappa(t)2b^2\mathbb{E}\left[\frac{X_t^2}{(b^2-X_t^2)^2}\right] - \frac{b^4}{W_{\text{e}}}\mathbb{E}\left[\frac{X_t^2}{(b^2-X_t^2)^3}\right] + \frac{b^2}{W_{\text{e}}} \mathbb{E}\left[\frac{b^2+3X_t^2}{(b^2-X_t^2)^3}\right].
\]
To complete our choice of four macroscopic state variables, we choose the first two moments popping up in the evolution equation of $\mathbb{E}[f_2(X_t)]$: $\mathbb{E}\left[\frac{X_t^2}{(b^2-X_t^2)^2}\right]$ and $\mathbb{E}\left[\frac{X_t^2}{(b^2-X_t^2)^3}\right]$.

\section{Derivation of first order convergence of approximate model~\eqref{eq:linearapproximate}} \label{app:B}
In this Appendix, we show that the $L_2$-error between the slow means of a general system of linear SDEs with an external period driving force, and the corresponding approximate macroscopic model for the slow component, decreases linearly to zero when the small-scale parameter $\varepsilon$ decreases to zero. For the derivation, we consider a more general form of the system in Section~\ref{subsec:periodic}
\begin{equation*}
\begin{aligned}
dX &= -\lambda(X+Y)dt + E\sin(at)dt + dW_x \\
dY &= \frac{1}{\varepsilon}(X-Y)dt + \frac{1}{\sqrt{\varepsilon}} dW_y,
\end{aligned}    
\end{equation*}
where the associated approximate macroscopic model is given by
\begin{equation*}
d\overset{\_}{X} = -2\lambda\overset{\_}{X}dt +E\sin(a t)dt+ dW_x.
\end{equation*}
One can show using elementary techniques from calculus that the mean vector $(\mu_X(t), \mu_Y(t))^T$ propagates as
\begin{equation} \label{eq:lineardrivensolution}
\begin{pmatrix} \mu_X(t) \\ \mu_Y(t) \end{pmatrix} = e^{t M} \begin{pmatrix} \mu_{X_0}-A \\ \mu_{Y_0}-C \end{pmatrix} + \begin{pmatrix}A \\ C \end{pmatrix} \cos(at) + \begin{pmatrix} B \\ D \end{pmatrix} \sin(at),
\end{equation}
where $M = \begin{pmatrix} -\lambda & -\lambda \\ \frac{1}{\varepsilon} & -\frac{1}{\varepsilon} \end{pmatrix}$ and $\begin{pmatrix}\mu_{X_0} \\ \mu_{Y_0} \end{pmatrix}$ is the initial condition to the equation. The constants $A,B,C$ and $D$ are the solution of the linear system 
\begin{equation*}
\begin{pmatrix} -a & \lambda & 0 & \lambda \\ \lambda & a & \lambda & 0 \\ 0 & -\frac{1}{\varepsilon} & -a & \frac{1}{\varepsilon} \\ -\frac{1}{\varepsilon} & 0 & \frac{1}{\varepsilon} & a \end{pmatrix} \begin{pmatrix} A \\ B \\ C \\ D \end{pmatrix} = \begin{pmatrix} E \\ 0 \\0 \\0 \end{pmatrix}.
\end{equation*}
Since we are only interested in the evolution of $\mu_X(t)$, we only need to know the explicit expressions for $A$ and $B$. Using the symbolic engine MUPAD \cite{fuchssteiner1996mupad}, the constants $A$ and $B$ are
\[
\begin{aligned}
A &= \frac{E\left(a(\lambda \varepsilon-1) - a^3\varepsilon^2\right)}{a^2(\lambda \varepsilon - 1)^2 + 4\lambda^2 + a^4\varepsilon^2} \\
B &= \frac{E\lambda (a^2\varepsilon^2+2)}{a^2(\lambda\varepsilon-1)^2 +4\lambda^2 + a^4\varepsilon^2}.
\end{aligned}
\]
Similarly, the evolution of the mean of the approximate macroscopic model is given by 
\begin{equation}
\overset{\_}{\mu}_X(t) = \left(\mu_{X_0}+ \frac{aE}{a^2+4\lambda^2}\right)e^{-2\lambda t}  - \frac{aE}{a^2+4\lambda^2}\cos(at) + \frac{2\lambda E}{a^2+4\lambda^2}\sin(at).
\end{equation}

We will first compute the difference between the means $\mu_X(t)$ and $\overset{\_}{\mu}_X(t)$ before giving an argument that the $L_2$-error decreases linearly with $\varepsilon$. The difference is means is given by
\begin{equation*}
\begin{aligned}
\mu_X(t) - \overset{\_}{\mu}_X(t) &= \left(\frac{E\left(a(\lambda \varepsilon-1) - a^3\varepsilon^2\right)}{a^2(\lambda \varepsilon - 1)^2 + 4\lambda^2 + a^4\varepsilon^2}+\frac{aE}{a^2+4\lambda^2}\right) \cos(at) \\
&+ \left(\frac{E\lambda (a^2\varepsilon^2+2)}{a^2(\lambda\varepsilon-1)^2 +4\lambda^2 + a^4\varepsilon^2} - \frac{2\lambda E}{a^2+4\lambda^2}\right) \sin(at) \\
&= a E \frac{(\lambda \varepsilon-1 - a^2\varepsilon^2)(a^2+4\lambda^2) + a^2(\lambda \varepsilon - 1)^2 + 4\lambda^2 + a^4\varepsilon^2}{(a^2(\lambda \varepsilon - 1)^2 + 4\lambda^2 + a^4\varepsilon^2)(a^2+4\lambda^2)} \cos(at) \\
&+ \lambda E \frac{(a^2\varepsilon^2+2)(a^2+4\lambda^2) - 2a^2(\lambda\varepsilon-1)^2 - 8\lambda^2 - 2a^4\varepsilon^2}{(a^2(\lambda \varepsilon - 1)^2 + 4\lambda^2 + a^4\varepsilon^2)(a^2+4\lambda^2)} \sin(at)
\end{aligned}
\end{equation*}
After some calculations, a lot of terms cancel in the expression above and the difference between the means reduces to
\begin{equation*}
\begin{aligned}
   \mu_X(t) - \overset{\_}{\mu}_X(t) &=  \varepsilon a E \frac{4\lambda^3 - \lambda^2 a - 3a^2 \lambda^2 \varepsilon}{(a^2(\lambda \varepsilon - 1)^2 + 4\lambda^2 + a^4\varepsilon^2)(a^2+4\lambda^2)} \cos(at) \\ 
   &+ \varepsilon \lambda E \frac{4a^2\lambda - a^4 \varepsilon + 2a^2 \lambda^2 \varepsilon}{(a^2(\lambda \varepsilon - 1)^2 + 4\lambda^2 + a^4\varepsilon^2)(a^2+4\lambda^2)} \sin(at).
  \end{aligned}
\end{equation*}
The dominant terms in both numerators linearly depend on $\varepsilon$, and the leading term in the common denominator is independent of $\varepsilon$. The difference in means hence decreases linearly with $\varepsilon$, and consequently the $L_2$-norm of the difference in means too, proving our statement in Section~\ref{subsec:periodic}.

\end{document}